\documentclass[12pt,leqno]{article}
\usepackage{amssymb}
\usepackage[mathscr]{eucal}
\usepackage{amsmath,amssymb,latexsym,theorem,bbm}
\newcommand{\DS}{\displaystyle}

\newcommand{\SC}{\scriptstyle}

\newcommand{\CC}{\mathsf{C}}
\newcommand{\DD}{\mathsf{D}}

\newcommand{\NN}{\mathbb{N}}

\newcommand{\RR}{\mathbb{R}}

\newcommand{\ZZ}{\mathbb{Z}}

\newcommand{\bA}{{\boldsymbol{A}}}
\newcommand{\tbA}{\widetilde{\bA}}
\newcommand{\bb}{{\boldsymbol{b}}}

\newcommand{\bd}{{\boldsymbol{d}}}

\newcommand{\bY}{{\boldsymbol{Y}}}

\newcommand{\bSigma}{{\boldsymbol{\Sigma}}}

\newcommand{\cA}{{\mathcal A}}
\newcommand{\cB}{{\mathcal B}}

\newcommand{\cF}{{\mathcal F}}

\newcommand{\cL}{{\mathcal L}}

\newcommand{\cN}{{\mathcal N}}

\newcommand{\cU}{{\mathcal U}}
\newcommand{\cV}{{\mathcal V}}

\newcommand{\cY}{{\mathcal Y}}
\newcommand{\cW}{{\mathcal W}}

\newcommand{\dd}{\mathrm{d}}

\newcommand{\slu}{{\SC\mathrm{lu}}}

\newcommand{\INARtwo}{\textup{INAR(2)}}
\newcommand{\INARp}{\textup{INAR($p$)}}

\newcommand{\EE}{\operatorname{E}}
\newcommand{\PP}{\operatorname{P}}

\newcommand{\halpha}{\widehat{\alpha}}
\newcommand{\hbeta}{\widehat{\beta}}

\newcommand{\tK}{\widetilde{K}}

\newcommand{\tsigma}{\widetilde{\sigma}}

\newcommand{\vare}{\varepsilon}
\newcommand{\mvb}{\mathversion{bold}}

\renewcommand{\mid}{\,|\,}

\renewcommand{\leq}{\leqslant}
\renewcommand{\geq}{\geqslant}

\newcommand{\stoch}{\stackrel{\PP}{\longrightarrow}}
\newcommand{\distr}{\stackrel{\cL}{\longrightarrow}}
\newcommand{\distre}{\stackrel{\cL}{=}}

\newcommand{\qmean}{\stackrel{L_2}{\longrightarrow}}

\newcommand{\lu}{\stackrel{\slu}{\longrightarrow}}
\newcommand{\as}{\stackrel{{\mathrm{a.s.}}}{\longrightarrow}}
\newcommand{\ase}{\stackrel{{\mathrm{a.s.}}}{=}}

\newcommand{\bbone}{\mathbbm{1}}

\newcommand{\nt}{{\lfloor nt\rfloor}}

\newcommand{\proofend}{\hfill\mbox{$\Box$}}

\numberwithin{equation}{section}

\theoremstyle{change} \theorembodyfont{\em}
\newtheorem{Lem}{Lemma.}[section]
\newtheorem{Thm}{Theorem.}[section]
\newtheorem{Pro}{Proposition.}[section]

\newtheorem{Def}{Definition.}[section]

\theorembodyfont{\rm}
\newtheorem{Rem}{Remark.}[section]

\begin{document}

\begin{center}
 {\bfseries\Large Asymptotic behavior of CLS estimator of autoregressive parameter
                   for nonprimitive unstable \INARtwo\ models} \\[5mm]

 {\sc\large M\'aty\'as $\text{Barczy}^{*,\diamond}$,
            \ M\'arton $\text{Isp\'any}^*$, \ Gyula $\text{Pap}^{\star}$}
\end{center}

\vskip0.2cm

\noindent * Faculty of Informatics, University of Debrecen,
            Pf.~12, H--4010 Debrecen, Hungary.

\noindent $\star$ Bolyai Institute, University of Szeged,
            Aradi v\'ertan\'uk tere 1, H--6720 Szeged, Hungary.

\noindent e--mails: barczy.matyas@inf.unideb.hu (M. Barczy),
                    ispany.marton@inf.unideb.hu (M. Isp\'any),
                    papgy@math.u-szeged.hu (G. Pap).

\noindent $\diamond$ Corresponding author.

%\vskip0.2cm

%\centerline{\sl March 06, 2009.}

\renewcommand{\thefootnote}{}
\vspace*{0.2cm}
\footnote{The authors have been supported by the Hungarian
 Portuguese Intergovernmental S \& T Cooperation Programme for 2008-2009
 under Grant No.\ PT-07/2007.
 M. Barczy and G. Pap have been partially supported by the Hungarian Scientific
 Research Fund under Grant No.\ OTKA T-079128.}

\vspace*{-10mm}

\begin{abstract}
In this paper the asymptotic behavior of conditional least squares estimators
 of the autoregressive parameter for nonprimitive unstable integer-valued autoregressive
 models of order 2 (\INARtwo) is described.
\end{abstract}

\section{Introduction and main results}\label{section_introduction_main_results}

Recently, there has been remarkable interest in integer-valued
 time series models (especially from statistical point of views) and
 a number of results are now available in specialized monographs  and review papers
 (e.g., Steutel and van Harn \cite{SteHar} and Wei{\ss} \cite{Wei1}).
Reasons to introduce discrete data models come from the need to account for
 the discrete nature of certain data sets, often counts of events, objects or
 individuals.

Among the most successful integer-valued time series models proposed in the
 literature we mention the INteger-valued AutoRegressive model of order $p$
 \ (\INARp). This model was first introduced by McKenzie \cite{McK}
 and Al-Osh and Alzaid \cite{AloAlz1} for the case $p=1$.
The INAR(1) model has been investigated by several authors.
The more general \INARp\ processes were first introduced by Al-Osh and Alzaid
 \cite{AloAlz2}.
In their setup the autocorrelation structure of the process corresponds to that of
 an ARMA($p,p-1$) process.
Another definition of an INAR(p) process was proposed independently
 by Du and Li \cite{DuLi} and by Gauthier and Latour \cite{GauLat} and Latour \cite{Lat2},
 and is different from that of Alzaid and Al-Osh \cite{AloAlz2}.
In Du and Li's setup the autocorrelation structure
 of an \INARp\ process is the same as that of an AR($p$) process. The setup of Du and Li
 \cite{DuLi} has been followed by most of the authors, and our approach will also be
 the same.
In Barczy et al. \cite{BarIspPap0} we investigated the asymptotic behavior of unstable INAR($p$)
 processes, i.e., when the characteristic polynomial has a unit root.
Under some natural assumptions we proved that the sequence of appropriately scaled random step
 functions formed from an unstable \INARp\ process converges weakly towards a squared Bessel process.
This limit process is a continuous branching process  also known as square-root process or
 Cox-Ingersoll-Ross process.

Parameter estimation for INAR models has a long history.
Franke and Seligmann \cite{FraSel} analyzed conditional maximum likelihood estimator
 of some parameters (including the autoregressive parameter)
 for stable INAR(1) models with Poisson innovations.
Du and Li \cite[Theorem 4.2]{DuLi} proved asymptotic normality of
 the conditional least squares (CLS) estimator of the autoregressive parameter for stable \INARp \ models,
 Br\"{a}nn\"{a}s and Hellstr\"{o}m \cite{BraHel} considered generalized method of moment estimation.
Silva and Oliveira \cite{SilOli} proposed a frequency domain based estimator of the autoregressive
 parameter for stable INAR(p) models with Poisson innovations.
Isp\' any et al. \cite{IspPapZui1} derived asymptotic inference for nearly unstable INAR(1) models
 which has been refined by Drost et al. \cite{DroAkkWer2} later.
Drost et al. \cite{DroAkkWer1} studied asymptotically efficient estimation of the parameters for
 stable INAR(p) models.

In this paper the asymptotic behavior of CLS estimators
 of the autoregressive parameter for so called nonprimitive unstable \INARtwo\ models is described,
 see our main results Theorem \ref{10main} and Theorem \ref{01main} later on.
In a forthcoming paper we will study asymptotic behavior of CLS estimators
 of the autoregressive parameter for primitive unstable \INARtwo\ models.
%There are also attempts to study the asymptotic behavior of CLS estimators for unstable \INARp \ models.

Concerning relevance and practical applications of unstable INAR models we note that
 empirical studies show importance of these kind of models.
Br\"{a}nn\"{a}s and Hellstr\"{o}m  \cite{BraHel} reported an INAR$(0.98)$ model for the number
 of private schools, Rudholm \cite{Rud} considered INAR$(0.98)$ \ and INAR$(0.99)$ \ models for
 the number of Swedish generic-pharmaceutical market.
Hellstr\"{o}m \cite{Hel} focused on the testing of unit root in INAR(1) models and provided
 small sample distributions for the Dickey-Fuller test statistic under the null hypothesis
 of unit root in an INAR(1) model with Poisson distributed innovations.
To our knowledge a unit root test for general INAR$(p)$ models is not known, and from
 this point of view studying unstable INAR$(p)$ models is an important preliminary task.

First we recall INAR(2) models.
Let \ $\ZZ_+$, $\NN$, \ $\RR$ \ and \ $\RR_+$ \ denote the set of non-negative integers,
 positive integers, real numbers and non-negative real numbers, respectively.
Every random variable will be defined on a fixed probability space
 \ $(\Omega,\cA,\PP)$.

\begin{Def}
Let \ $(\vare_k)_{k\in\NN}$ \ be an independent and identically distributed
 (i.i.d.) sequence of non-negative integer-valued random variables, and let
 \ $(\alpha,\beta)\in [0,1]^2$.
\ An \INARtwo\ time series model with autoregressive parameter \ $(\alpha,\beta)\in[0,1]^2$
 \ and innovations \ $(\vare_k)_{k\in\NN}$ \ is a stochastic process
 \ $(X_k)_{k \geq -1}$ \ given by
 \begin{align}\label{INAR2}
   X_k = \sum_{j=1}^{X_{k-1}} \xi_{k,j}
         + \sum_{j=1}^{X_{k-2}} \eta_{k,j} + \vare_k , \qquad k \in \NN ,
 \end{align}
 where for all \ $k\in\NN$, \ $(\xi_{k,j})_{j\in\NN}$ \ and \ $(\eta_{k,j})_{j\in\NN}$
 \ are sequences of i.i.d.\ Bernoulli random variables with mean \ $\alpha$
 \ and \ $\beta$, \ respectively \ such that these sequences are mutually
 independent and independent of the sequence \ $(\vare_k)_{k\in\NN}$, \ and
 \ $X_0$, \ $X_{-1}$ \ are non-negative integer-valued random variables
 independent of the sequences \ $(\xi_{k,j})_{j\in\NN}$,
 \ $(\eta_{k,j})_{j\in\NN}$, \ $k\in\NN$, \ and \ $(\vare_k)_{k\in\NN}$.
\end{Def}

The \INARtwo\ model \eqref{INAR2} can be written in another way using
 the binomial thinning operator \ $\alpha\,\circ$
 \ (due to Steutel and van Harn~\cite{SteHar}) which we recall now.
Let \ $X$ \ be a non-negative integer-valued random variable.
Let \ $(\xi_j)_{j\in\NN}$ \ be a sequence of i.i.d.\ Bernoulli random variables
 with mean \ $\alpha\in[0,1]$.
\ We assume that the sequence \ $(\xi_j)_{j\in\NN}$ \ is independent of \ $X$.
\ The non-negative integer-valued random variable \ $\alpha\,\circ X$
 \ is defined by
 \[
   \alpha\circ X
     :=\begin{cases}
        \sum\limits_{j=1}^X\xi_j, & \quad \text{if \ $X>0$},\\[2mm]
         0, & \quad \text{if \ $X=0$}.
       \end{cases}
 \]
The sequence \ $(\xi_j)_{j\in\NN}$ \ is called a counting sequence.
The \INARtwo\ model \eqref{INAR2} takes the form
 \[
    X_k = \alpha \circ X_{k-1} + \beta \circ X_{k-2} + \vare_k ,
    \qquad k \in \NN .
 \]
Note that the above form of the INAR(2) model is quite analogous with a usual AR($2$) process
 (another slight link between them is the similarity of some conditional expectations,
 see \eqref{10seged1}).
This definition of the INAR(2) process was proposed independently
 by Du and Li \cite{DuLi} and by Gauthier and Latour \cite{GauLat} and Latour \cite{Lat2},
 and is different from that of Alzaid and Al-Osh \cite{AloAlz2}, which assumes that the
 conditional distribution of the vector
 \ $(\alpha\circ X_t, \beta\circ X_t)$ \ given \ $X_t = x_t$
 \ is multinomial with parameters \ $(\alpha, \beta, x_t)$ \ and
 is independent of the past history of the process.
The two different formulations imply different second-order structure for the processes:
 under the first approach, the INAR(2) has the same second-order structure as an AR($2$)
 process, whereas under the second one, it has the same one as an ARMA$(2,1)$
 process.

Based on the asymptotic behavior of \ $\EE(X_k)$ \ as \ $k\to\infty$ \
 described in Proposition 2.2 in Barczy et al.~\cite{BarIspPap0},
 we distinguish three types of INAR(2) models.
The asymptotic behavior of \ $\EE(X_k)$ \ as \ $k\to\infty$ \ is determined by the
 spectral radius \ $\varrho(\bA)$ \ of the matrix
 \[
   \bA
    :=\begin{bmatrix}
      \alpha & \beta \\
      1 & 0 \\
      \end{bmatrix},
 \]
 i.e., by the maximum of the modulus of the eigenvalues of \ $\bA$.
\ The case \ $\varrho(\bA)<1$, \ when \ $\EE(X_k)$ \ converges to a finite limit as \ $k\to\infty$,
 \ is called \emph{stable} or \emph{asymptotically stationary},
 whereas the cases \ $\varrho(\bA)=1$, \ when \ $\EE(X_k)$ \ tends linearly to \ $\infty$,
 \ and \ $\varrho(\bA)>1$, \ when \ $\EE(X_k)$ \ converges to \ $\infty$ \ with an exponential rate,
 are called \emph{unstable} and \emph{explosive}, respectively.
Clearly, \ $\varrho(\bA)<1$, \ $\varrho(\bA)=1$ \ and \ $\varrho(\bA)>1$ \ are equivalent with
 \ $\alpha + \beta < 1$, \ $\alpha + \beta = 1$ \ and \ $\alpha + \beta > 1$, \ respectively,
 see Barczy et al.~\cite[Proposition 2.1]{BarIspPap0}.

%If \ $\beta=0$, \ then \ $\alpha$ \ and \ $0$ \ are the eigenvalues of \ $\bA$, \ and hence
% the above mentioned equivalence holds also for the case \ $\beta=0$.

An INAR(2) process with autoregressive parameter \ $(\alpha,\beta)$ \ such that \ $\alpha > 0$ \ and
 \ $\beta > 0$ \ is called \emph{primitive}, otherwise, i.e., if \ $\alpha = 0$ \ or \ $\beta = 0$,
 \ it is called
 \emph{nonprimitive} (see Barczy et al.~\cite[Definition 2.2]{BarIspPap0}).
If \ $\alpha>0$ \ and \ $\beta=0$, \ then \ $(X_n)_{n\geq -1}$ \ is an INAR(1) process
 with autoregressive parameter \ $\alpha$.
\ If \ $\alpha=0$ \ and \ $\beta>0$, \ then \ $(X_n)_{n\geq -1}$ \ takes the form
 \[
    X_n=\beta \circ X_{n-2}+\vare_n,\quad n\in\NN,
 \]
 and hence the subsequences \ $(X_{2n-j})_{n\geq 0}$, $j=0,1$, \ form independent primitive
 INAR(1) processes with autoregressive parameter \ $\beta$ \ such that \ $X_{-j}=0$.

For the sake of simplicity we consider a zero start \INARtwo\ process, that is
 we suppose \ $X_0 = X_{-1} = 0$.
\ The general case of nonzero initial values may be handled in a similar way, but we renounce
 to consider it.

In the sequel we always assume that \ $\EE(\vare_1^2) < \infty$.
\ Let us denote the mean and variance of \ $\vare_1$ \ by \ $\mu_\vare$ \ and
 \ $\sigma_\vare^2$, \ respectively.
In all what follows we suppose that \ $\mu_\vare>0$, \ otherwise \ $X_k=0$ \ for all
 \ $k\in\NN$.

Next we formulate our main results considering the two nonprimitive unstable cases separately.
For all \ $n\in\NN$, \ a CLS estimator \ $(\halpha_n, \hbeta_n)$ \ of the
 autoregressive parameter \ $(\alpha,\beta)\in[0,1]^2$ \ based on a sample \ $X_1,\ldots,X_n$ \
 will be denoted by \ $(\halpha_n({\bf X}_n), \hbeta_n({\bf X}_n))$.
In Section \ref{section_CLS_estimators} we present a result about the existence and uniqueness of
 \ $(\halpha_n({\bf X}_n), \hbeta_n({\bf X}_n))$, \ see Proposition \ref{Pro1}.

\begin{Thm}\label{10main}
Let \ $(X_k)_{k \geq -1}$ \ be a nonprimitive \INARtwo\ process with autoregressive parameter
 \ $(1,0)$ \ (hence it is unstable).
Suppose that \ $X_0 = X_{-1} = 0$, \ $\EE(\vare_1^4) < \infty$ \ and
 \ $\mu_\vare > 0$.
\ Then
 \[
   \begin{bmatrix}
    \sqrt{n}(\halpha_n({\bf X}_n) - 1) \\
    \sqrt{n}\hbeta_n({\bf X}_n)
   \end{bmatrix}
   \distr \frac{2 \sigma_\vare}{\sqrt{\mu_\vare^2 + 4 \sigma_\vare^2}}
                       Z \begin{bmatrix} -1 \\ 1 \end{bmatrix}
   \qquad \text{as \ $n\to\infty$,}
 \]
 where \ $Z$ \ is a standard normally distributed random variable and
 \ $\distr$ \ denotes convergence in distribution.
Hence the limit distribution is a centered normal distribution with covariance matrix
 \[
    \frac{4 \sigma_\vare^2}{\mu_\vare^2 + 4 \sigma_\vare^2}
       \begin{bmatrix}
         1 & -1 \\
         -1 & 1 \\
       \end{bmatrix}.
 \]
\end{Thm}

The proof of Theorem \ref{10main} can be found in Section \ref{Nonprimitive_10_section}.

\begin{Rem}
We note that a fourth order moment condition on the innovation distribution in Theorem \ref{10main}
 is supposed (i.e., we suppose \ $\EE(\vare_1^4)<\infty$), \ which is used for checking the so called
 conditional Lindeberg condition of a martingale central limit theorem (see the proof of Theorem \ref{10main}).
However it is important to remark that this condition is a technical one,
 we suspect that Theorem \ref{10main} remains true under second order moment condition
 on the innovation distribution, but we renounce to consider it.
\proofend
\end{Rem}

\begin{Thm}\label{01main}
Let \ $(X_k)_{k \geq -1}$ \ be a nonprimitive \INARtwo\ process with autoregressive parameter
  \ $(0,1)$ \ (hence it is unstable).
Suppose that \ $X_0 = X_{-1} = 0$, \ $\EE(\vare_1^2) < \infty$ \ and
 \ $\mu_\vare > 0$.
\ Then
 \[
   \begin{bmatrix}
    n \halpha_n({\bf X}_n) \\
    n (\hbeta_n({\bf X}_n) - 1)
   \end{bmatrix}
   \distr
   \frac{\int_0^1 \cW_t \, \dd \cW_t}{\int_0^1 (\cW_t)^2 \, \dd t}
   \begin{bmatrix} -1 \\ 1 \end{bmatrix}
   \qquad \text{as \ $n\to\infty$,}
 \]
 where \ $(\cW_t)_{t \in \RR_+}$ \ is a standard Wiener process.
\end{Thm}

The proof of Theorem \ref{01main} can be found in Section \ref{Nonprimitive_01_section}.

\begin{Rem}
We recall that the distribution of \ $\int_0^1 \cW_t \, \dd \cW_t/\int_0^1 (\cW_t)^2 \, \dd t$
 \ is the same as the limit distribution of the Dickey-Fuller statistics, see, e.g.,
 the Ph.D. Thesis of Bobkoski \cite{Bob}, or (7.14) and Theorem 9.5.1 in Tanaka \cite{Tan}.
\proofend
\end{Rem}

\begin{Rem}
We note that in both nonprimitive unstable cases the limit distributions are concentrated on 
 the same line \ $\{(x,y)\in\RR^2 : x+y=0\}$.
\ However, these limit distributions are different. 
In the unstable case \ $(1,0)$ \ we have a centred normal limit distribution and 
 the difference of the CLS estimator \ $(\halpha_n({\bf X}_n),\hbeta_n({\bf X}_n))$ 
 \ and \ $(1,0)$ \ has to be normalized by \ $\sqrt{n}$. 
\ In the unstable case \ $(0,1)$ \ we have a different limit distribution
 (described in Theorem \ref{01main}) and we have to normalize by \ $n$ \ instead of \ $\sqrt{n}$.
\proofend
\end{Rem}

The rest of the paper is organized as follows.
In Section \ref{section_CLS_estimators} we study the CLS estimator
 of the autoregressive parameter \ $(\alpha,\beta)$ \ of nonprimitive unstable \INARtwo\ models.
Section \ref{Nonprimitive_10_section} and Section \ref{Nonprimitive_01_section} are devoted to the proofs
 considering the two nonprimitive unstable cases \ $(\alpha,\beta)=(1,0)$ \ and
  \ $(\alpha,\beta)=(0,1)$ \ separately.

\section{CLS estimators}\label{section_CLS_estimators}

For all \ $k \in \ZZ_+$, \ let us denote by \ $\cF_k$ \ the \ $\sigma$--algebra
 generated by the random variables \ $X_0,X_1,\ldots,X_k$.
\ (Note that \ $\cF_0 = \{ \Omega, \emptyset \}$, \ since \ $X_0 = 0$.)
By \eqref{INAR2},
 \begin{align}\label{10seged1}
   \EE(X_k\mid\cF_{k-1})
   = \alpha X_{k-1} + \beta X_{k-2} + \mu_\vare , \qquad k \in \NN .
 \end{align}
Let us introduce the sequence
 \begin{equation}\label{Mk}
  M_k := X_k - \EE(X_k \mid \cF_{k-1})
       = X_k - \alpha X_{k-1} - \beta X_{k-2} - \mu_\vare ,
  \qquad k \in \NN ,
 \end{equation}
 of martingale differences with respect to the filtration
 \ $(\cF_k)_{k \in \ZZ_+}$.
\ The process \ $(X_k)_{k \geq -1}$ \ satisfies the recursion
 \begin{equation}\label{regr}
  X_k = \alpha X_{k-1} + \beta X_{k-2} + M_k + \mu_\vare ,
  \qquad k \in \NN .
 \end{equation}
For all \ $n\in\NN$, \ a CLS estimator \ $(\halpha_n, \hbeta_n)$ \ of the
 autoregressive parameter \ $(\alpha,\beta)\in[0,1]^2$ \ based on a sample \ $X_1,\ldots,X_n$ \
 can be obtained by minimizing the sum of squares
 \begin{align}\label{SS1}
   \sum_{k=1}^n \big( X_k - \EE(X_k \mid \cF_{k-1}) \big)^2
   = \sum_{k=1}^n (X_k - \alpha X_{k-1} - \beta X_{k-2} - \mu_\vare)^2
 \end{align}
 with respect to \ $(\alpha,\beta)$ \ over \ $\RR^2$.
For all \ $n\in\NN$ \ and \ $x_1,\ldots,x_n\in\RR$, \ let us put
 \ ${\bf x}_n := (x_1,\ldots,x_n)$.
\ Motivated by \eqref{SS1}, for all \ $n\in\NN$, \ we define the function
 \ $Q_n:\RR^n\times\RR^2\to\RR$ \ by
 \[
    Q_n({\bf x}_n;\alpha',\beta') := \sum_{k=1}^n (x_k - \alpha' x_{k-1} - \beta' x_{k-2} - \mu_\vare)^2
 \]
 for all \ $\alpha',\beta'\in\RR$ \ and \ ${\bf x}_n\in\RR^n$ \ with \ $x_{-1}:=x_0:=0$.
\ By definition, for all \ $n\in\NN$, \ a CLS estimator of the autoregressive parameter
 \ $(\alpha,\beta)\in[0,1]^2$ \ is a measurable function
 \ $(\halpha_n, \hbeta_n):\RR^n\to\RR^2$ \ such that
 \[
     Q_n({\bf x}_n;\halpha_n({\bf x}_n),\hbeta_n({\bf x}_n))
        = \inf_{(\alpha',\beta')\in\RR^2} Q_n({\bf x}_n;\alpha',\beta')
        \qquad \forall\;\;  {\bf x}_n\in\RR^n.
 \]
For all \ $n\in\NN$ \ and \ $\omega\in\Omega$, \ let us put
 \begin{align*}
   {\bf X}_n(\omega) := (X_1(\omega),\ldots,X_n(\omega)),
    \qquad
   {\bf X}_n := (X_1,\ldots,X_n).
 \end{align*}

Next we give the explicit form of the CLS estimators \ $(\halpha_n, \hbeta_n)$, $n\in\NN$.

\begin{Lem}\label{Lem_CLSE}
Any measurable function \ $(\halpha_n, \hbeta_n):\RR^n\to\RR^2$ \ for which
 \begin{align}\label{seged1}
   \begin{bmatrix}
     \halpha_n({\bf x}_n)  \\
     \hbeta_n({\bf x}_n)
    \end{bmatrix}
     =  \begin{bmatrix}
           \sum_{k=1}^n x_{k-1}^2 &  \sum_{k=1}^n x_{k-1}x_{k-2} \\
           \sum_{k=1}^n x_{k-1}x_{k-2} &  \sum_{k=1}^n x_{k-2}^2
         \end{bmatrix}^{-1}
             \begin{bmatrix}
               \sum_{k=1}^n( x_k - \mu_\vare )x_{k-1} \\
               \sum_{k=1}^n( x_k - \mu_\vare )x_{k-2}
             \end{bmatrix},
 \end{align}
 if \ $\sum_{k=1}^n x_{k-2}^2>0$, \ and
 \[
     \halpha_n({\bf x}_n) = \frac{x_n-\mu_\vare}{x_{n-1}},
 \]
 if \ $\sum_{k=1}^n x_{k-2}^2=0$ \ and \ $x_{n-1}\ne0$,
 \ is a CLS estimator of the autoregressive parameter \ $(\alpha,\beta)\in[0,1]^2$.
\end{Lem}

We note that \ $(\halpha_n, \hbeta_n)$ \ is not defined uniquely on the set
 \ $\{{\bf x}_n\in\RR^n : \sum_{k=1}^n x_{k-2}^2 = 0 \}$.

\noindent
\textbf{Proof of Lemma \ref{Lem_CLSE}.}
First we note that for all \ $({\bf x}_n;\alpha',\beta')\in\RR^n\times\RR^2$,
\begin{align*}
   &\frac{\partial Q_n}{\partial \alpha'}({\bf x}_n;\alpha',\beta')
     = -2\sum_{k=1}^n \big(x_k-\alpha' x_{k-1} - \beta'x_{k-2}-\mu_\vare\big)x_{k-1}, \\
   &\frac{\partial Q_n}{\partial \beta'}({\bf x}_n;\alpha',\beta')
     = -2\sum_{k=1}^n \big(x_k-\alpha' x_{k-1} - \beta'x_{k-2}-\mu_\vare\big)x_{k-2}, \\
   &\frac{\partial^2 Q_n}{\partial (\alpha')^2}({\bf x}_n;\alpha',\beta')
     = 2\sum_{k=1}^n x_{k-1}^2,
     \qquad
    \frac{\partial^2 Q_n}{\partial (\beta')^2}({\bf x}_n;\alpha',\beta')
      = 2\sum_{k=1}^n x_{k-2}^2, \\
   &\frac{\partial^2 Q_n}{\partial \alpha'\partial \beta'}({\bf x}_n;\alpha',\beta')
     = 2\sum_{k=1}^n x_{k-1}x_{k-2}.
 \end{align*}

Now let us suppose that \ $\sum_{k=1}^n x_{k-2}^2>0$.
\ It is enough to show that the function
 \[
    \RR^2\ni(\alpha',\beta')\mapsto Q_n({\bf x}_n;\alpha',\beta')
 \]
 is strictly convex and that \eqref{seged1} is the unique solution of the system of equations
 \begin{align}\label{SS1_der_eq2}
   \frac{\partial Q_n}{\partial \alpha'}({\bf x}_n;\alpha',\beta') = 0,
    \qquad
   \frac{\partial Q_n}{\partial \beta'}({\bf x}_n;\alpha',\beta') = 0.
 \end{align}
In proving strict convexity of the function in question, it is enough to check that
 the \ $(2\times 2)$ \ Hessian matrix
 \[
      \begin{bmatrix}
          \frac{\partial^2 Q_n}{\partial (\alpha')^2} & \frac{\partial^2 Q_n}{\partial \beta'\partial\alpha'} \\
          \frac{\partial^2 Q_n}{\partial \alpha'\partial\beta'} & \frac{\partial^2 Q_n}{\partial (\beta')^2} \\
        \end{bmatrix}
        ({\bf x}_n;\alpha',\beta')
        =
          \begin{bmatrix}
           2\sum_{k=1}^n x_{k-1}^2 &  2\sum_{k=1}^n x_{k-1}x_{k-2} \\
           2\sum_{k=1}^n x_{k-1}x_{k-2} &  2\sum_{k=1}^n x_{k-2}^2
         \end{bmatrix}
 \]
 is (strictly) positive definite, see, e.g., Berkovitz \cite[Theorem 3.3, Chapter III]{Ber}.
Since \ $\sum_{k=1}^n x_{k-2}^2>0$, \ there exists some \ $i\in\{1,\ldots,n-2\}$
 \ such that \ $x_i\ne0$ \ and hence there does not exist a constant \ $c\in\RR$ \ such that
  \ $(x_0,x_1,\ldots,x_{n-1}) = c(x_{-1},x_0,\ldots,x_{n-2})$.
\ Then \ $(x_0,x_1,\ldots,x_{n-1})$ \ and \ $(x_{-1},x_0,\ldots,x_{n-2})$ \ are linearly independent,
 and, by Cauchy and Schwarz's inequality, we get
 \[
     \sum_{k=1}^n x_{k-1}^2 \sum_{k=1}^n x_{k-2}^2
        > \left(\sum_{k=1}^n x_{k-1}x_{k-2}\right)^2.
 \]
Hence the above \ $(2\times 2)$ \ Hessian matrix has positive leading principal minors and
 then it is positive definite.
An easy calculation shows that \eqref{seged1} satisfies \eqref{SS1_der_eq2}.

Now let us suppose that \ $\sum_{k=1}^n x_{k-2}^2=0$ \ and \ $x_{n-1}\ne0$.
\ Then
 \begin{align}\label{01seged19}
   Q_n({\bf x}_n;\alpha',\beta')
       = (x_n-\alpha'x_{n-1}-\mu_\vare)^2 + (x_{n-1}-\mu_\vare)^2+ (n-2)\mu_\vare^2
       \qquad \forall\;\, (\alpha',\beta')\in\RR^2,
 \end{align}
 and for all \ $(\alpha',\beta')\in\RR^2$,
 \begin{align*}
  &\frac{\partial Q_n}{\partial \alpha'}({\bf x}_n;\alpha',\beta')
     = -2(x_n-\alpha'x_{n-1}-\mu_\vare)x_{n-1},\\
  &\frac{\partial Q_n}{\partial \beta'}({\bf x}_n;\alpha',\beta')
     =0.
 \end{align*}
An easy calculation shows that for any function \ $\hbeta_n:\RR^n\to\RR$, \
 \[
   \begin{bmatrix}
     \frac{x_n-\mu_\vare}{x_{n-1}} \\[1mm]
     \hbeta_n({\bf x}_n) \\
   \end{bmatrix}
 \]
 is a solution of \eqref{SS1_der_eq2}.
By \eqref{01seged19}, \ $Q_n$ \ as a function of \ $\alpha'$ \ is a polynomial of order 2,
 and hence \ $(x_n-\mu_\vare)/x_{n-1}$ \ is a global minimum of \ $Q_n$ \ (as function of \ $\alpha'$).

Finally, let us suppose that \ $\sum_{k=1}^n x_{k-2}^2=0$ \ and \ $x_{n-1}=0$.
\ Then
  \[
   Q_n({\bf x}_n;\alpha',\beta')
     = (x_n-\mu_\vare)^2 + (n-1)\mu_\vare^2,
     \qquad \forall\;\, (\alpha',\beta')\in\RR^2,
  \]
 which yields the statement.
\proofend

In the sequel by the expression {\sl `a property holds
 asymptotically as \ $n\to\infty$ \ with probability one'} we mean that there
 exists an event \ $S\in\cA$ \ such that \ $\PP(S)=1$ \ and for all
 \ $\omega\in S$ \ there exists an \ $n(\omega)\in\NN$ \ such that the property
 in question holds for all \ $n\geq n(\omega)$.
\ Next we present a result about the existence and uniqueness of
 \ $(\halpha_n({\bf X}_n), \hbeta_n({\bf X}_n))$.

\begin{Pro}\label{Pro1}
Let \ $(X_k)_{k \geq -1}$ \ be a nonprimitive \INARtwo\ process with autoregressive parameter
  \ $(1,0)$ \ or \ $(0,1)$.
\ Suppose that \ $X_0=X_{-1}=0$, \ $\EE(\vare_1^2)<\infty$ \ and \ $\mu_\vare>0$.
\ Then the following statements hold
 asymptotically as \ $n\to\infty$ \ with probability one:
 \ $\sum_{k=1}^n X_{k-2}^2>0$ \ and hence there exists a unique CLS estimator
 \ $(\halpha_n({\bf X}_n), \hbeta_n({\bf X}_n))$ \ having the form
  \begin{equation}\label{CLSE}
  \begin{bmatrix} \halpha_n({\bf X}_n) \\ \hbeta_n({\bf X}_n) \end{bmatrix} = \bA_n^{-1} \bb_n ,
 \end{equation}
 where
 \[
   \bA_n := \sum_{k=1}^n
             \begin{bmatrix}
              X_{k-1}^2 & X_{k-1}X_{k-2} \\
              X_{k-1}X_{k-2} & X_{k-2}^2
             \end{bmatrix} , \qquad
   \bb_n := \sum_{k=1}^n
             \begin{bmatrix}
              ( X_k - \mu_\vare )X_{k-1} \\
              ( X_k - \mu_\vare )X_{k-2}
             \end{bmatrix} .
 \]
%If  \ $(\alpha,\beta)=(0,1)$, \ then
% \[
%   \lim_{n\to\infty}\PP\left(\sum_{k=1}^n X_{k-2}^2>0\right)=1,
% \]
% and hence the probability of the existence of a unique CLS estimator
% \ $(\halpha_n({\bf X}_n), \hbeta_n({\bf X}_n))$ \ converges to 1 as \ $n\to\infty$,
% \ and this CLS estimator has the form given in \eqref{CLSE}.
\end{Pro}

\noindent{\bf Proof.}
First we consider the case of \ $(1,0)$.
\ In this case equation \eqref{INAR2} has the form \ $X_k=X_{k-1}+\vare_k$, $k\in\NN$, \ and hence
 \ $X_n=\sum_{k=1}^n\vare_k$, $n\in\NN$.
\ By the strong law of large numbers we have
 \begin{align}
  n^{-1} X_n = n^{-1} \sum_{k=1}^n \vare_k \as \mu_\vare ,
  \label{SLLNXn}
 \end{align}
 and hence
 \[
   n^{-2} X_n^2 \as \mu_\vare^2 ,
 \]
 where \ $\as$ \ denotes almost sure convergence.
Then \ $X_n/n^3\as0$ \ and \ $X_n^2/n^3\as0$, \ and hence, by Toeplitz theorem, we conclude
 \begin{equation}\label{10SumX2}
  n^{-3} \sum_{k=1}^n X_{k-2}^2 \as \frac{1}{3} \mu_\vare^2 .
 \end{equation}
Since \ $\mu_\vare>0$, \ by \eqref{10SumX2}, we get \ $\sum_{k=1}^n X_{k-2}^2>0$ \
 holds asymptotically as \ $n\to\infty$ \ with probability one and Lemma \ref{Lem_CLSE}
 yields that there exists a unique CLS estimator
 \ $(\halpha_n({\bf X}_n), \hbeta_n({\bf X}_n))$ \ having the form \eqref{CLSE}
 asymptotically as \ $n\to\infty$ \ with probability one.

Next we consider the case of \ $(0,1)$.
\ In this case equation \eqref{INAR2} has the form \ $X_k=X_{k-2}+\vare_k$, $k\in\NN$, \ and hence
 \ $X_{2n}=\sum_{k=1}^n\vare_{2k}$, $n\in\ZZ_+$, \ and \ $X_{2n-1}=\sum_{k=1}^n\vare_{2k-1}$, $n\in\ZZ_+$.
\ By the strong law of large numbers, we have
 \[
  n^{-1} X_{2n} \as \mu_\vare, \quad \text{as \ $n\to\infty$,}
  \qquad \text{and}\qquad
  n^{-1} X_{2n-1} \as \mu_\vare \quad \text{as \ $n\to\infty$,}
 \]
 which yield that
 \[
  n^{-1} X_{n} \as \frac{1}{2}\mu_\vare \quad \text{as \ $n\to\infty$.}
 \]
Using Toeplitz theorem, as in the case of \ $(1,0)$, \ we get
 \[
  n^{-3} \sum_{k=1}^n X_{k-2}^2 \as \frac{1}{12} \mu_\vare^2 .
 \]
One can finish the proof as in the case of \ $(1,0)$.
\proofend

In Section \ref{Nonprimitive_10_section} and Section \ref{Nonprimitive_01_section}
 we will usually write \ $(\halpha_n,\hbeta_n)$ \ instead of
 \ $(\halpha_n({\bf X}_n), \hbeta_n({\bf X}_n))$.

\section{Proofs for the nonprimitive unstable case \textbf{\mvb $(1,0)$}}
\label{Nonprimitive_10_section}

In the case of \ $(\alpha, \beta) =  (1,0)$, \ equation \eqref{INAR2} has the form \ $X_k = X_{k-1} + \vare_k$,
 \ $k \in \NN$, \ hence in fact, we have a random walk
 \ $X_k = \vare_1 + \cdots + \vare_k$, \ $k \in \NN$, \ with positive drift \ $\mu_\vare$, \ since
 \ $\EE(X_k) = \mu_\vare k$, \ $k \in \NN$.

%The next lemma describes the asymptotic behavior of the sequence \ $(\det({\bf A}_n))_{n\in\NN}$.

%\begin{Lem}\label{Lem10}
%Let \ $(X_k)_{k \geq -1}$ \ be a nonprimitive \INARtwo\ process with coefficients
% \ $(\alpha, \beta) =  (1,0)$ \ (hence it is unstable).
%Suppose that \ $X_0 = X_{-1} = 0$, \ $\EE(\vare_1^2) < \infty$ \ and \ $\mu_\vare>0$.
%\ Then
% \begin{equation}\label{10det}
%  n^{-4} \det(\bA_n)
%  \as \frac{1}{12} \mu_\vare^2 \left( 4\sigma_\vare^2 + \mu_\vare^2 \right)
%    \qquad \text{as \ $n\to\infty$,}
% \end{equation}
% where \ $\as$ \ denotes almost sure convergence.
%\end{Lem}

Next we present an auxiliary lemma which will be used in the proof of Theorem \ref{10main}.

\begin{Lem}\label{Lem10_seged}
Let \ $\xi_n$, $\eta_n$, $n\in\NN$ \ and \ $\xi$ \ be random variables such that
 \ $\xi_n\distr\xi$ \ as \ $n\to\infty$ \ and
 \ $\lim_{n\to\infty} \PP(\xi_n=\eta_n)=1$.
\ Then \ $\eta_n\distr\xi$ \ as \ $n\to\infty$.
\end{Lem}

\noindent{\bf Proof.}
We give three proofs. \ Let \ $x\in\RR$ \ be a continuity point of the distribution function of \ $\xi$.
 \ Then for all \ $n\in\NN$,
 \begin{align*}
   \PP(\eta_n<x)
     & = \PP(\eta_n<x,\xi_n=\eta_n) + \PP(\eta_n<x,\xi_n\ne\eta_n) \\
     & = \PP(\xi_n<x,\xi_n=\eta_n) + \PP(\eta_n<x,\xi_n\ne\eta_n).
 \end{align*}
Since \ $\PP(\eta_n<x,\xi_n\ne\eta_n)\leq \PP(\xi_n\ne\eta_n)$, \ we have
 \ $\lim_{n\to\infty} \PP(\eta_n<x,\xi_n\ne\eta_n) = 0$ \ and
 \begin{align*}
   \lim_{n\to\infty} \PP(\xi_n<x,\xi_n=\eta_n)
    & = \lim_{n\to\infty} \left( \PP(\xi_n<x) - \PP(\xi_n<x,\xi_n\ne\eta_n)\right)\\
    & = \lim_{n\to\infty} \PP(\xi_n<x)= \PP(\xi<x).
 \end{align*}
Hence \ $\lim_{n\to\infty} \PP(\eta_n<x)= \PP(\xi<x)$.

Our second proof sounds as follows.
For all \ $\vare>0$, \ we have
 \begin{align*}
   \PP(\vert\eta_n-\xi_n\vert\geq \vare)
    & = \PP(\vert\eta_n-\xi_n\vert\geq \vare,\eta_n=\xi_n)
       + \PP(\vert\eta_n-\xi_n\vert\geq \vare,\eta_n\ne\xi_n)\\
    & = \PP(\vert\eta_n-\xi_n\vert\geq \vare,\eta_n\ne\xi_n).
 \end{align*}
Since \ $\lim_{n\to\infty} \PP(\xi_n=\eta_n)=1$, \ we have
 \[
    \lim_{n\to\infty} \PP(\vert\eta_n-\xi_n\vert\geq \vare,\eta_n\ne\xi_n) = 0,
 \]
 and hence \ $\lim_{n\to\infty}\PP(\vert\eta_n-\xi_n\vert\geq \vare) = 0$ \ $\forall$ $\vare>0$, \ i.e.,
 \ $\eta_n-\xi_n$ \ converges in probability to \ $0$ \ as \ $n\to\infty$.
\ Then Slutsky's lemma yields the assertion.

Our third proof sounds as follows.
For all \ $\vare>0$, \ we have \ $\PP(\vert\eta_n-\xi_n\vert\geq\vare)\leq\PP(\eta_n\ne\xi_n)$, $n\in\NN$,
 \ which yields that \ $\eta_n-\xi_n$ \ converges in probability to \ $0$ \ as \ $n\to\infty$.
\ Then Slutksky's lemma yields the assertion.
\proofend

\vskip0.5cm

\noindent
\textbf{Proof of Theorem \ref{10main}.}
By Proposition \ref{Pro1},
 \[
   \begin{bmatrix}
     \halpha_n - \alpha \\
     \hbeta_n - \beta
    \end{bmatrix} = \bA_n^{-1} \bd_n
 \]
 holds asymptotically as \ $n\to\infty$ \ with probability one, where
 \begin{align}\label{10def_dn}
   \bd_n := \sum_{k=1}^n
             \begin{bmatrix}  M_k X_{k-1} \\  M_k X_{k-2} \end{bmatrix},
             \qquad n\in\NN.
 \end{align}
We can write
 \[
   \bA_n^{-1} \bd_n = \frac{1}{\det(\bA_n)} \tbA_n \bd_n
 \]
 asymptotically as \ $n\to\infty$ \ with probability one,
 where \ $\tbA_n$ \ denotes the adjoint of \ $\bA_n$ \ given by
 \[
   \tbA_n
   :=\sum_{k=1}^n
      \begin{bmatrix}
       X_{k-2}^2 & - X_{k-1}X_{k-2} \\
       - X_{k-1}X_{k-2} & X_{k-1}^2
      \end{bmatrix} .
 \]

Next we study the asymptotic behavior of the sequence \ $(\det(\bA_n))_{n \in \NN}$.
Namely, we show that
 \begin{equation}\label{10det}
   n^{-4} \det(\bA_n)
   \as \frac{1}{12} \mu_\vare^2 \left( 4\sigma_\vare^2 + \mu_\vare^2 \right)
    \qquad \text{as \ $n\to\infty$.}
 \end{equation}
We note that for deriving \eqref{10det} we need only second order moment condition on the innovation
 distribution (i.e., \ $\EE(\vare_1^2)<\infty$), \ the fourth order moment condition
 \ $\EE(\vare_1^4)<\infty$ \ will be used in the description of the asymptotic behavior
 of the sequence \ $(\tbA_n \bd_n)_{n\in\NN}$.
\ We have
 \begin{align}
  \begin{split}
  \det(\bA_n) &= \sum_{k=1}^n X_{k-1}^2 \sum_{k=1}^n X_{k-2}^2
                 - \left( \sum_{k=1}^n X_{k-1} X_{k-2} \right)^2 \\
             &= \sum_{k=1}^n (X_{k-2} + \vare_{k-1})^2 \sum_{k=1}^n X_{k-2}^2
                 - \left( \sum_{k=1}^n (X_{k-2} + \vare_{k-1}) X_{k-2} \right)^2 \\
             &= \sum_{k=1}^n X_{k-2}^2 \sum_{k=1}^n \vare_{k-1}^2
                 - \left( \sum_{k=1}^n X_{k-2} \, \vare_{k-1} \right)^2 ,
  \end{split}
  \label{detAn}
 \end{align}
 where \ $\vare_0:=0$.
\ By the strong law of large numbers we have
 \begin{equation}\label{10Sumvare2}
  n^{-1} \sum_{k=1}^n \vare_{k-1}^2
  \as \EE(\vare_1^2) = \sigma_\vare^2 + \mu_\vare^2 .
 \end{equation}
Moreover,
 \[
   \sum_{k=1}^n X_{k-2} \, \vare_{k-1}
   = \sum_{k=1}^n \vare_{k-1} \sum_{i=1}^{k-2} \vare_i
   = \sum_{1 \leq i < j \leq n-1} \vare_i \vare_j
   = \frac{1}{2}
     \left( \left( \sum_{k = 1}^n \vare_{k-1} \right)^2
            - \sum_{k = 1}^n \vare_{k-1}^2 \right) ,
 \]
 and hence, by \eqref{SLLNXn} and \eqref{10Sumvare2},
 \begin{equation}\label{10SumXvare}
  n^{-2} \sum_{k=1}^n X_{k-2} \, \vare_{k-1}
  \as \frac{1}{2} \mu_\vare^2 .
 \end{equation}
By \eqref{detAn}, \eqref{10SumX2}, \eqref{10Sumvare2} and \eqref{10SumXvare}, we
 deduce \eqref{10det}.

Now we study the asymptotic behavior of the sequence
 \ $(\tbA_n \bd_n)_{n \in \NN}$.
\ First note that \ $M_k = X_k - X_{k-1} - \mu_\vare = \vare_k - \mu_\vare$,
 \ $k \in \NN$, \ since \ $\alpha=1$ \ and \ $\beta=0$.
\ We have
 \begin{align*}
  \tbA_n \bd_n
  &= \begin{bmatrix}
      \sum\limits_{k=1}^n X_{k-2}^2
       & - \sum\limits_{k=1}^n ( X_{k-2} + \vare_{k-1} ) X_{k-2} \\[3mm]
      - \sum\limits_{k=1}^n ( X_{k-2} + \vare_{k-1} ) X_{k-2}
       & \sum\limits_{k=1}^n ( X_{k-2} + \vare_{k-1} )^2
     \end{bmatrix}
     \begin{bmatrix}
      \sum\limits_{k=1}^n (\vare_k - \mu_\vare) ( X_{k-2} + \vare_{k-1} ) \\[3mm]
      \sum\limits_{k=1}^n (\vare_k - \mu_\vare) X_{k-2}
     \end{bmatrix} \\
  &= e_n^{(1)}
     \begin{bmatrix}
      1 \\ -1
     \end{bmatrix}
     +
     e_n^{(2)}
      \begin{bmatrix}
      0 \\ -1
     \end{bmatrix} ,
 \end{align*}
 where
 \begin{align*}
   & e_n^{(1)} := \sum\limits_{k=1}^n X_{k-2}^2
               \sum\limits_{k=1}^n (\vare_k - \mu_\vare) \vare_{k-1}
               - \sum\limits_{k=1}^n  \vare_{k-1} X_{k-2}
                 \sum\limits_{k=1}^n (\vare_k - \mu_\vare) X_{k-2},\quad n\in\NN,\\[1mm]
   & e_n^{(2)} := \sum\limits_{k=1}^n \vare_{k-1} X_{k-2}
                  \sum\limits_{k=1}^n (\vare_k - \mu_\vare) \vare_{k-1}
                  - \sum\limits_{k=1}^n \vare_{k-1}^2
                   \sum\limits_{k=1}^n (\vare_k - \mu_\vare) X_{k-2},\quad n\in\NN.
 \end{align*}

The aim of the following discussion is to apply multidimensional martingale central limit
 theorem (see, e.g., Jacod and Shiryaev \cite[Chapter VIII, Theorem 3.33]{JSh}) for the sequences
 \ $(\bY_{n,k}, \, \cF_k)_{k \in \NN}$, \ $n \in \NN$, \ of square-integrable
 martingale differences, where
 \[
   \bY_{n,k} := \begin{bmatrix}
               n^{-3/2} (\vare_k - \mu_\vare) X_{k-2} \\
               n^{-1/2} (\vare_k - \mu_\vare) \vare_{k-1}
              \end{bmatrix} , \qquad n, k \in \NN ,
 \]
 where \ $\vare_0=0$.
\ Using that the \ $\sigma$-algebra generated by \ $\vare_1,\ldots,\vare_k$ \ equals
 \ $\cF_k$ \ for all \ $k\in\NN$, \ we get \ $\EE\big(\bY_{n,k}\mid \cF_{k-1}\big)={\bf 0}\in\RR^2$
 \ and
 \[
   \EE \big( \bY_{n,k} \bY_{n,k}^\top \mid \cF_{k-1} \big)
   = \sigma_\vare^2
     \begin{bmatrix}
      n^{-3} X_{k-2}^2 & n^{-2} X_{k-2} \, \vare_{k-1} \\
      n^{-2} X_{k-2} \, \vare_{k-1} & n^{-1} \vare_{k-1}^2
     \end{bmatrix} ,
  \qquad n,k\in\NN.
 \]
Hence by \eqref{10SumX2}, \eqref{10Sumvare2} and \eqref{10SumXvare} we have the
 asymptotic covariance matrices
 \[
   \sum_{k=1}^\nt
    \EE \big( \bY_{n,k} \bY_{n,k}^\top \mid \cF_{k-1} \big)
   \as \sigma_\vare^2
       \begin{bmatrix}
        \frac{t^3}{3} \mu_\vare^2 & \frac{t^2}{2} \mu_\vare^2 \\
        \frac{t^2}{2} \mu_\vare^2 & t(\sigma_\vare^2 + \mu_\vare^2)
       \end{bmatrix}
       =: \bSigma(t) , \qquad t\in\RR_+ ,
 \]
 where \ $\lfloor x\rfloor$ \ denotes the integer part of a real number \ $x\in\RR$.
\ The conditional Lindeberg condition
 \[
   \sum_{k=1}^\nt
    \EE \big( \|\bY_{n,k}\|^2 \bbone_{\{ \|\bY_{n,k}\| > \theta \}} \mid \cF_{k-1} \big)
   \stoch 0
 \]
 is satisfied for all \ $t\in\RR_+$ \ and \ $\theta > 0$, \ where \ $\stoch$ \ denotes
 convergence in probability.
Indeed, using that \ $\EE(\vare_1^4)<\infty$,
 \begin{align*}
  \sum_{k=1}^\nt
    \EE \big( \|\bY_{n,k}\|^2 \bbone_{\{ \|\bY_{n,k}\| > \theta \}} \mid \cF_{k-1} \big)
  &\leq \frac{1}{\theta^2} \sum_{k=1}^\nt
        \EE \big( \|\bY_{n,k}\|^4 \mid \cF_{k-1} \big) \\
  &\leq \frac{2}{\theta^2} \sum_{k=1}^\nt
        \EE \big( n^{-6} (\vare_k - \mu_\vare)^4 X_{k-2}^4
                  + n^{-2} (\vare_k - \mu_\vare)^4 \vare_{k-1}^4 \mid \cF_{k-1}
            \big) \\
  &= \frac{2 \EE\big[(\vare_1 - \mu_\vare)^4\big]}{\theta^2} \sum_{k=1}^\nt
     \big( n^{-6} X_{k-2}^4 + n^{-2} \vare_{k-1}^4 \big)
   \stoch 0 ,
 \end{align*}
 where the last step follows by \ $\EE(X_k^4)\leq k^4\EE(\vare_1^4)$, $k\in\NN$.
\ Indeed, by power mean inequality
 \begin{align*}
  \frac{X_k}{k} = \frac{1}{k}\sum_{i=1}^k\vare_i
                \leq  \left(\frac{1}{k}\sum_{i=1}^k\vare_i^4\right)^{1/4},
                \qquad k\in\NN,
 \end{align*}
 and hence
  \[
     k^{-4}\EE(X_k^4) \leq \frac{1}{k}\sum_{i=1}^k\EE(\vare_i^4) = \EE(\vare_1^4),
       \qquad k\in\NN.
  \]
Thus we obtain
 \[
   \sum_{k=1}^n\bY_{n,k}
    =\sum_{k=1}^n
      \begin{bmatrix}
       n^{-3/2} (\vare_k - \mu_\vare) X_{k-2} \\
       n^{-1/2} (\vare_k - \mu_\vare) \vare_{k-1}
      \end{bmatrix}
      \distr
      \cN\left(\begin{bmatrix}
              0 \\
              0 \\
            \end{bmatrix},
      \bSigma(1)\right) .
 \]
By \eqref{10SumX2}, \eqref{10SumXvare} and Slutsky's lemma, we obtain
 \begin{align*}
  n^{-7/2}e_n^{(1)}
  = \sum_{k=1}^n
      \begin{bmatrix}
       - n^{-2} X_{k-2} \, \vare_{k-1} \\
       n^{-3} X_{k-2}^2
      \end{bmatrix}^\top
     \sum_{k=1}^n
      \begin{bmatrix}
       n^{-3/2} (\vare_k - \mu_\vare) X_{k-2} \\
       n^{-1/2} (\vare_k - \mu_\vare) \vare_{k-1}
      \end{bmatrix}
   \distr \cN(0,\sigma^2) ,
 \end{align*}
 where
 \[
   \sigma^2 := \begin{bmatrix}
               - \frac{1}{2} \mu_\vare^2 \\
               \frac{1}{3} \mu_\vare^2
              \end{bmatrix}^\top
              \bSigma(1)
              \begin{bmatrix}
               - \frac{1}{2} \mu_\vare^2 \\
               \frac{1}{3} \mu_\vare^2
              \end{bmatrix}
            = \frac{1}{36} \mu_\vare^4 \sigma_\vare^2
              (\mu_\vare^2 + 4 \sigma_\vare^2) .
 \]
In a similar way, by \eqref{10Sumvare2}, \eqref{10SumXvare} and Slutsky's lemma,
 \begin{align*}
  n^{-5/2} e_n^{(2)}
  = \sum_{k=1}^n
      \begin{bmatrix}
       - n^{-1} \vare_{k-1}^2 \\
       n^{-2} X_{k-2} \vare_{k-1}
      \end{bmatrix}^\top
     \sum_{k=1}^n
      \begin{bmatrix}
       n^{-3/2} (\vare_k - \mu_\vare) X_{k-2} \\
       n^{-1/2} (\vare_k - \mu_\vare) \vare_{k-1}
      \end{bmatrix}
   \distr \cN(0,\tsigma^2) ,
 \end{align*}
 where
 \[
   \tsigma^2 := \begin{bmatrix}
                 - ( \mu_\vare^2 + \sigma_\vare^2 ) \\
                 \frac{1}{2} \mu_\vare^2
                \end{bmatrix}^\top
                \bSigma(1)
                \begin{bmatrix}
                 - ( \mu_\vare^2 + \sigma_\vare^2 ) \\
                 \frac{1}{2} \mu_\vare^2
                \end{bmatrix}
              = \frac{1}{12} \mu_\vare^2 \sigma_\vare^2
                (\mu_\vare^2 + \sigma_\vare^2) (\mu_\vare^2 + 4 \sigma_\vare^2) .
 \]
Then, by Slutsky's lemma, \ $ n^{-7/2} e_n^{(2)}\distr 0$ \ as \ $n\to\infty$, \ which also yields
 that  \ $ n^{-7/2} e_n^{(2)}\stoch 0$ \ as \ $n\to\infty$.
\ Consequently, again by Slutsky's lemma,
 \[
   n^{-7/2} \tbA_n \bd_n
      \distr \sigma Z \begin{bmatrix}
                       1 \\ -1
                      \end{bmatrix},
  \]
 where \ $Z$ \ is a standard normally distributed random variable.
Using part (v) of Theorem 2.7 in van der Vaart \cite{Vaa}, \eqref{10det} yields that
 \begin{align*}
   \left(n^{-4}\det(\bA_n), n^{-7/2} \tbA_n \bd_n \right)
     \distr
   \left( \frac{1}{12} \mu_\vare^2 (4\sigma_\vare^2 + \mu_\vare^2),
           \sigma Z \begin{bmatrix}
                       1 \\ -1
                      \end{bmatrix}
  \right)
  \qquad \text{ as \ $n\to\infty$.}
 \end{align*}
Let us introduce the function \ $g:\RR\times\RR^2\to\RR^2$,
 \begin{align}\label{01seged12}
  g\left(x,\begin{bmatrix}
        y \\
        z \\
      \end{bmatrix}\right)
   :=\begin{cases}
            \begin{bmatrix}
              y/x \\
              z/x \\
            \end{bmatrix}, & \text{if \ $x\ne 0$,}\\[5mm]
            \begin{bmatrix}
              0 \\
              0 \\
            \end{bmatrix}, & \text{if \ $x=0$.}
            \end{cases}
 \end{align}
Since \ $g$ \ is continuous on \ $(\RR\setminus\{0\})\times\RR^2$ \ and
 \[
   \PP\left(
      \left( \frac{1}{12} \mu_\vare^2 (4\sigma_\vare^2 + \mu_\vare^2),
           \sigma Z \begin{bmatrix}
                       1 \\ -1
                      \end{bmatrix}
    \right)  \in (\RR\setminus\{0\})\times\RR^2
   \right)=1,
 \]
 the continuous mapping theorem (see, e.g., Theorem 2.3 in van der Vaart \cite{Vaa})
 yields that
 \[
   g\left(n^{-4}\det(\bA_n), n^{-7/2} \tbA_n \bd_n \right)
     \distr
         \frac{12\sigma}{\mu_\vare^2(\mu_\vare^2 + 4\sigma_\vare^2)}
          Z \begin{bmatrix}
                       1 \\ -1
            \end{bmatrix}
       =\frac{2\sigma_\vare}{\sqrt{\mu_\vare^2 + 4\sigma_\vare^2}}
          Z \begin{bmatrix}
                       1 \\ -1
            \end{bmatrix}
 \]
 as \ $n\to\infty$.
\ By Proposition \ref{Pro1}, we have
 \begin{align*}
  \PP\left(
      \sqrt{n}\begin{bmatrix}
                \halpha_n-\alpha \\
                \hbeta_n -\beta \\
              \end{bmatrix}
      = g\left(n^{-4}\det(\bA_n), n^{-7/2} \tbA_n \bd_n \right)
     \right)
   \geq \PP\left(\sum_{k=1}^n X_{k-2}^2 > 0\right)
   \to 1
   \qquad \text{as \ $n\to \infty$.}
 \end{align*}
Then Lemma \ref{Lem10_seged} yields the assertion.
\proofend

%\noindent{\bf\large Proofs for nonprimitive unstable case \ $(\alpha, \beta) =  (0,1)$}

\section{Proofs for the nonprimitive unstable case \textbf{\mvb $(0,1)$}}
\label{Nonprimitive_01_section}

The structure of the proof is the same as the proof of Theorem \ref{10main}
 (nonprimitive unstable case \ $(1,0))$.
\ Namely, based on the decomposition
 \begin{align}\label{01seged20}
    \begin{bmatrix}
     \halpha_n - \alpha \\
     \hbeta_n - \beta
    \end{bmatrix}
    = \frac{1}{\det(\bA_n)} \tbA_n \bd_n,
 \end{align}
 which holds asymptotically as \ $n\to\infty$ \ with probability one (see Proposition \ref{Pro1}),
 first we will study the asymptotic behavior of the sequence \ $(\det(\bA_n))_{n\in\NN}$ \ and then
 the asymptotic behavior of the sequence \ $(\tbA_n \bd_n)_{n\in\NN}$.
\ The main differences from the proof of Theorem \ref{10main} are that the reference to the
 strong law of large numbers and Toeplitz theorem in the case of \ $(\det(\bA_n))_{n\in\NN}$,
 \ and reference to the multidimensional martingale central limit theorem in the case
 \ $(\tbA_n \bd_n)_{n\in\NN}$ \ should be replaced and completed here by, for example, (asymptotic)
 expansions separating the expectations ('leading terms') of the entries of \ $\bA_n$ \ and
 the coordinates of \ $\tbA_n \bd_n$, \ respectively.
In the case of \ $(1,0)$ \ it was proved that \ $n^{-4}\det(\bA_n)$ \ converges
 almost surely to a positive non-random limit (see \eqref{10det}) and hence, by the decomposition
 \eqref{01seged20}, to prove convergence in distribution of the appropriately normalized
 sequence
 \[
    \begin{bmatrix}
     \sqrt{n}(\halpha_n - 1) \\
     \sqrt{n}\hbeta_n
    \end{bmatrix},\quad n\in\NN,
 \]
 it was enough to prove convergence in distribution of the appropriately normalized sequence
 \ $n^{-7/2}(\tbA_n \bd_n)_{n\in\NN}$.
\ In contrast to the case \ $(1,0)$ \ it will turn out that \ $n^{-6}\det(\bA_n)$ \
 converges almost surely to 0 (see \eqref{A}) in the case of \ $(0,1)$,
 \ and hence the method used for the case \ $(1,0)$ \ can not be carried out
 in the case of \ $(0,1)$.
\ However, we can prove that \ $n^{-5}\det(\bA_n)$ \ converges in distribution to a positive random limit
 (see Lemma \ref{Lem01}) and \ $n^{-4}\tbA_n \bd_n$ \ converges also in distribution
 (see the proof of Theorem \ref{10main}).
To be able to use the decomposition \eqref{01seged20}, we need to establish joint convergence
 in distribution of \ $n^{-5}\det(\bA_n)$ \ and \ $n^{-4}\tbA_n \bd_n$.
\ For this reason we will derive (asymptotic) expansions for \ $\det(\bA_n)$, \ $\tbA_n$ \ and \ $\bd_n$,
 \ respectively, such that these expansions will consist of the same 'building blocks'.
These 'building blocks' are listed in Lemma \ref{Lem01_UV3} and their joint convergence in
 distribution is also proved which yields that \ $n^{-5}\det(\bA_n)$ \ and \ $n^{-4}\tbA_n \bd_n$ \
 also converge jointly in distribution.
To prove Lemma \ref{Lem01_UV3}, using multidimensional martingale central theorem,
 we will verify that
  \[
    \begin{bmatrix}
      n^{-1/2}(X_{2n} - \EE(X_{2n})) \\
       n^{-1/2}(X_{2n-1} - \EE(X_{2n-1})) \\
    \end{bmatrix}
  \]
 converges in distribution as \ $n\to\infty$ \ (see Lemma \ref{Lem01_UV})
 and then an appropriate version of the continuous mapping theorem will be used.

First we recall two versions of the continuous mapping theorem for
 \ $\RR^d$-valued stochastic processes with c\`adl\`ag paths.

A function \ $f : \RR_+ \to \RR^d$ \ is called \emph{c\`adl\`ag} if it is right
 continuous with left limits.
\ Let \ $\DD(\RR_+, \RR^d)$ \ and \ $\CC(\RR_+, \RR^d)$ \ denote the space of all
 \ $\RR^d$-valued c\`adl\`ag and continuous functions on \ $\RR_+$, \ respectively.
Let \ $\cB(\DD(\RR_+, \RR^d))$ \ denote the Borel $\sigma$-field in \ $\DD(\RR_+, \RR^d)$
 \ for the metric defined in Jacod and Shiryaev \cite[Chapter VI, (1.26)]{JSh}
 (with this metric \ $\DD(\RR_+, \RR^d)$ \ is a complete and separable metric space
  and the topology induced by this metric is the so-called Skorokhod topology).
For \ $\RR^d$-valued stochastic processes \ $(\cY_t)_{t\in\RR_+}$ \ and \ $(\cY^n_t)_{t\in\RR_+}$,
 $n\in\NN$, \ with c\`adl\`ag paths we write \ $\cY^n \distr \cY$ \ if the
 distribution of \ $\cY^n$ \ on the space \ $(\DD(\RR_+, \RR^d),\cB(\DD(\RR_+, \RR^d)))$ \ converges weakly to
 the distribution of \ $\cY$ \ on the space \ $(\DD(\RR_+, \RR^d),\cB(\DD(\RR_+, \RR^d)))$ \ as \ $n\to\infty$.
\ Concerning the notation \ $\distr$ \ we note that if \ $\xi_n$, $n\in\NN$, \ and \ $\xi$ \ are
 random elements with values in a metric space \ $(E,d)$, \ then we also denote by \ $\xi_n\distr \xi$
 \  the weak convergence of the distributions of \ $\xi_n$ \ on the space \ $(E,\cB(E))$ \ towards
 the distribution of \ $\xi$ \ on the space \ $(E,\cB(E))$ \ as \ $n\to\infty$, \ where \ $\cB(E)$ \
 denotes the Borel \ $\sigma$-algebra on \ $E$ \ induced by the given metric \ $d$.

The following version of continuous mapping theorem can be found for example
 in Kallenberg \cite[Theorem 3.27]{Kal}.

\begin{Lem}\label{Lem_Kallenberg}
Let \ $(S,d_S)$ \ and \ $(T,d_T)$ \ be metric spaces and \ $(\xi_n)_{n\in\NN}$, \ $\xi$ \ be random elements
 with values in \ $S$ \ such that \ $\xi_n\distr\xi$ \ as \ $n\to\infty$.
\ Let \ $f:S\to T$ \ and \ $f_n:S\to T$, $n\in\NN$, \ be measurable mappings and
 \ $C\in\cB(S)$ \ such that \ $\PP(\xi\in C)=1$ \ and
 \ $\lim_{n\to\infty} d_T(f_n(s_n),f(s))=0$ \ if \ $\lim_{n\to\infty} d_S(s_n,s)=0$ \ and \ $s\in C$.
\ Then \ $f_n(\xi_n)\distr f(\xi)$ \ as \ $n\to\infty$.
\end{Lem}

For the case \ $S:=\DD(\RR_+, \RR^d)$ \ and \ $T:=\RR^q$, \ where \ $d$, $q\in\NN$ \ we formulate
 a consequence of Lemma \ref{Lem_Kallenberg}.

For a function \ $f \in \DD(\RR_+, \RR^d)$ \ and for a sequence \ $(f_n)_{n\in\NN}$
 \ in \ $\DD(\RR_+, \RR^d)$, \ we write \ $f_n \lu f$ \ if \ $(f_n)_{n\in\NN}$
 \ converges to \ $f$ \ locally uniformly, i.e., if
 \ $\sup_{t\in[0,T]} \Vert f_n(t) - f(t)\Vert \to 0$ \ as \ $n \to \infty$
 \ for all \ $T > 0$.
\ For measurable mappings \ $\Phi : \DD(\RR_+, \RR^d) \to \RR^q$ \ and
 \ $\Phi_n : \DD(\RR_+, \RR^d) \to \RR^q$, \ $n \in \NN$, \ we will
 denote by \ $C_{\Phi,(\Phi_n)_{n\in\NN}}$ \ the set of all functions
 \ $f \in \CC(\RR_+, \RR^d)$ \ such that \ $\Phi_n(f_n) \longrightarrow \Phi(f)$ \ whenever
 \ $f_n \lu f$ \ with \ $f_n \in \DD(\RR_+, \RR^d)$, \ $n \in \NN$.

\begin{Lem}\label{Conv2Funct}
Let \ $(\cU_t)_{t\in\RR_+}$ \ and \ $(\cU^n_t)_{t\in\RR_+}$, $n\in\NN$, \ be
 \ $\RR^d$-valued stochastic processes with c\`adl\`ag paths such that \ $\cU^n \distr \cU$
 \ as \ $n\to\infty$.
\ Let \ $\Phi : \DD(\RR_+, \RR^d) \to \RR^q$ \ and
 \ $\Phi_n : \DD(\RR_+, \RR^d) \to \RR^q$, $n \in \NN$, \ be
 measurable mappings such that there exists \ $C\subset C_{\Phi,(\Phi_n)_{n\in\NN}}$ \ with
 \ $C\in\cB(\DD(\RR_+, \RR^d))$ \ and \ $\PP(\cU \in C) = 1$.
\ Then \ $\Phi_n(\cU^n) \distr \Phi(\cU)$ \ as \ $n\to\infty$.
\end{Lem}

\noindent
\textbf{Proof.}
First we recall that for all \ $g\in C(\RR_+, \RR^d)$, \ $g_n\in \DD(\RR_+, \RR^d)$, $n\in\NN$, \ the sequence
 \ $(g_n)_{n\in\NN}$ \ converges to \ $g$ \ in the Skorokhod topology of
 \ $\DD(\RR_+, \RR^d)$ \ if and only if it converges to \ $g$ \ locally uniformly
 (see, e.g., Jacod and Shiryaev \cite[Chapter VI., Proposition 1.17. (b)]{JSh}), i.e.,
 with the notation \ $\stackrel{S_d}{\longrightarrow}$ \ for convergence in the Skorokhod topology
 of \ $D(\RR_+,\RR^d)$, \ $g_n\lu g$ \ if and only if \ $g_n\stackrel{S_d}{\longrightarrow} g$.
\ Hence
 \begin{align*}
  C_{\Phi,(\Phi_n)_{n\in\NN}}
            &= \left\{f\in C(\RR_+,\RR^d) : \Phi_n(f_n) \longrightarrow \Phi(f), \,\forall\, f_n \lu f,\,
                   f_n\in \DD(\RR_+, \RR^d), n\in\NN \right\}\\[1mm]
            &= \left\{ f\in C(\RR_+,\RR^d) : \Phi_n(f_n) \longrightarrow \Phi(f),
                                         \,\forall \, f_n \stackrel{S_d}{\longrightarrow} f, \,
                                         f_n\in \DD(\RR_+, \RR^d), n\in\NN \right\}.
 \end{align*}
Then Lemma \ref{Lem_Kallenberg} with the special choices
 \ $S:=D(\RR_+,\RR^d)$, \ $T:=\RR^q$, \ $\xi_n:=(\cU^n_t)_{t\in\RR_+}$,
 \ \ $\xi:=(\cU_t)_{t\in\RR_+}$, \ $f_n:=\Phi_n$, $n\in\NN$, \ and \ $f:=\Phi$ \
 yields the assertion.
\proofend

We also remark that a slightly different proof of Lemma \ref{Conv2Funct} can be found
 in Isp\'any and Pap \cite[Lemma 3.1]{IspPap}.

In the case of \ $(\alpha, \beta) =  (0,1)$, \ equation \eqref{INAR2} has the form \ $X_k = X_{k-2} + \vare_k$,
 \ $k \in \NN$, \ hence in fact, now we have two independent random walks
 \begin{align*}
  U_k &:= X_{2k} = \sum_{j=1}^k \vare_{2j}, \qquad k \in \ZZ_+ , \\
  V_k &:= X_{2k-1} = \sum_{j=1}^k \vare_{2j-1}, \qquad k \in \ZZ_+ ,
 \end{align*}
 with positive drifts \ $\mu_\vare$, \ since \ $\EE(U_k) = \mu_\vare k$,
 \ $k \in \ZZ_+$, \ and \ $\EE(V_k) = \mu_\vare k$, \ $k \in \ZZ_+$, \ respectively.
Let us introduce the random step functions
 \[
   \cU_t^n := U_\nt, \qquad \cV_t^n := V_\nt, \qquad t \in\RR_+, \qquad n \in \NN.
 \]

In what follows we present several lemmas which will be used later on.

\begin{Lem}\label{Lem01_UV}
Let \ $(X_k)_{k \geq -1}$ \ be a nonprimitive \INARtwo\ process with autoregressive parameter
  \ $(0,1)$ \ (hence it is unstable).
Suppose that \ $X_0 = X_{-1} = 0$, \ $\EE(\vare_1^2) < \infty$ and \ $\mu_\vare>0$.
\ Then
  \begin{align}\label{01seged1}
  \begin{bmatrix}
   n^{-1/2} \big(\cU^n - \EE(\cU^n) \big) \\
   n^{-1/2} \big(\cV^n - \EE(\cV^n) \big)
  \end{bmatrix}
  \distr
  \begin{bmatrix}
   \sigma_\vare \cW^{(1)} \\
   \sigma_\vare \cW^{(2)}
  \end{bmatrix},
 \end{align}
 where \ $(\cW^{(1)}_t)_{t\in\RR_+}$ \ and \  $(\cW^{(2)}_t)_{t\in\RR_+}$ \ are
 independent standard Wiener processes.
Further, for all \ $\delta>1/2$,
 \begin{align}\label{01seged9}
     \frac{U_n-\EE(U_n)}{n^\delta} \stoch 0,
       \qquad
     \frac{V_n-\EE(V_n)}{n^\delta} \stoch 0.
 \end{align}
\end{Lem}

\noindent{\bf Proof.}
We show that the multidimensional martingale central limit theorem (see, e.g., Jacod and Shiryaev
 \cite[Chapter VIII, Theorem 3.33]{JSh}) implies \eqref{01seged1}.
Indeed, with the notation
 \[
   \bY_{n,k} := \begin{bmatrix}
               n^{-1/2} (\vare_{2k} - \mu_\vare) \\
               n^{-1/2} (\vare_{2k-1} - \mu_\vare)
              \end{bmatrix} , \qquad n,k \in \NN ,
 \]
 we have \ $(\bY_{n,k},\cF_{2k})_{k\in\NN}$, $n\in\NN$, \ are sequences of square-integrable
 martingale differences such that \ $\EE\big(\bY_{n,k}\mid \cF_{2(k-1)}\big)={\bf 0}\in\RR^2$
 \ and
 \[
   \EE \big( \bY_{n,k} \bY_{n,k}^\top \mid \cF_{2(k-1)} \big)
   = \sigma_\vare^2
     n^{-1}
      I_2,
  \qquad n,k\in\NN,
 \]
 where \ $I_2$ \ denotes the $2\times 2$ identity matrix.
Then the asymptotic covariance matrices
 \[
   \sum_{k=1}^\nt
     \EE \big( \bY_{n,k} \bY_{n,k}^\top \mid \cF_{2(k-1)} \big)
      \as \sigma_\vare^2 t I_2,
      \qquad t\in\RR_+.
 \]
The conditional Lindeberg condition
  \begin{align}\label{01seged4}
   \sum_{k=1}^\nt
    \EE \big( \|\bY_{n,k}\|^2 \bbone_{\{ \|\bY_{n,k}\| > \theta \}} \mid \cF_{2(k-1)} \big)
   \stoch 0
  \end{align}
 is satisfied for all \ $t\in\RR_+$ \ and \ $\theta > 0$.
\ Indeed, we have
 \begin{align*}
  \sum_{k=1}^\nt
    & \EE \big( \|\bY_{n,k}\|^2 \bbone_{\{ \|\bY_{n,k}\| > \theta \}} \big) \\
    & = \frac{1}{n} \sum_{k=1}^\nt
        \EE\Big[\left( (\vare_{2k}-\mu_\vare)^2 + (\vare_{2k-1}-\mu_\vare)^2 \right)
                         \bbone_{\{ (\vare_{2k}-\mu_\vare)^2 + (\vare_{2k-1}-\mu_\vare)^2 > n\theta^2 \}}
            \Big] \\
    & = \frac{\nt}{n}
        \EE\Big[\left( (\vare_2-\mu_\vare)^2 + (\vare_1-\mu_\vare)^2 \right)
                         \bbone_{\{ (\vare_2-\mu_\vare)^2 + (\vare_1-\mu_\vare)^2 > n\theta^2 \}}
            \Big]
     \to 0,
   \end{align*}
 by dominated convergence theorem.
This yields that the convergence in \eqref{01seged4} holds in fact in $L_1$-sense.
Thus we obtain \eqref{01seged1}.
\proofend

\begin{Lem}\label{Lem01_UV2}
Let \ $d, p \in \NN$ \ and let \ $K : [0,1] \times \RR^d \to \RR^p$ \ be
 a function such that for all \ $R > 0$ \ there exists
 \ $C_R > 0$ such that
 \begin{equation}\label{01seged5}
  \| K(t, x) - K(s, y) \| \leq C_R \left( | t - s | + \| x - y \| \right)
 \end{equation}
 for all \ $s, t \in [0, 1]$ \ and \ $x, y \in \RR^d$ \ with \ $\| x \| \leq R$
 \ and \ $\| y \| \leq R$.
\ Moreover, let us define the mappings
 \ $\Phi, \Phi_n : \DD(\RR_+, \RR^d) \to \RR^{d + p}$, \ $n \in \NN$, \ by
 \begin{align*}
  \Phi_n(f)
  &:= \left( f(1),
             \frac{1}{n}
             \sum_{k=1}^n
              K\left( \frac{k}{n}, f\left( \frac{k}{n} \right) \right) \right),
  \\
  \Phi(f)
  &:= \left( f(1), \int_0^1 K( t, f(t) ) \, \dd t \right)
 \end{align*}
 for all \ $f \in \DD(\RR_+, \RR^d)$.
\ Then the mappings \ $\Phi_n$, $n \in \NN$, \ and  \ $\Phi$ \ are measurable, and
 \ $C_{\Phi,(\Phi_n)_{n \in \NN}} = \CC(\RR_+, \RR^d) \in \cB(\DD(\RR_+, \RR^d))$.
\end{Lem}

\noindent{\bf Proof.}
For an arbitrary Borel set \ $B \in \cB(\RR^{d + p})$ \ we have
 \[
   \Phi_n^{-1}(B) = \pi_{\frac{1}{n}, \frac{2}{n}, \dots, 1}^{-1}(\tK_n^{-1}(B)) ,
   \qquad n\in\NN,
 \]
 where for all \ $n\in\NN$ \ the mapping \ $\tK_n : (\RR^d)^n \to \RR^{d + p}$ \ is defined by
 \[
   \tK_n(x_1, \dots, x_n)
   := \left( x_n,
             \frac{1}{n} \sum_{k=1}^n K\left( \frac{k}{n}, x_k \right) \right) ,
   \qquad x_1, \dots, x_n \in \RR^d ,
 \]
 and the natural projections
 \ $\pi_{t_1, t_2, \dots, t_n} : \DD(\RR_+, \RR^d) \to (\RR^d)^n$,
 \ $t_1, t_2, \dots, t_n \in \RR_+$, \ are given by
 \ $\pi_{t_1, t_2, \dots, t_n}(f) := (f(t_1), f(t_2), \dots, f(t_n))$,
 \ $f \in \DD(\RR_+, \RR^d)$, \ $t_1, t_2, \dots, t_n \in \RR_+$.
 \ Since \ $K$ \ is continuous, \ $\tK_n$ \ is also continuous, and hence
 \ $\tK_n^{-1}(B)\in\cB((\RR^d)^n)$.
\ It is known that \ $\pi_{t_1, t_2, \dots, t_n}$, \ $t_1, t_2, \dots, t_n \in \RR_+$,
 \ are measurable mappings
 (see, e.g., Billingsley \cite[Theorem 16.6 (ii)]{Bil} or Ethier and Kurtz \cite[Proposition 3.7.1]{EthKur}),
 and hence \ $\Phi_n=\tK_n\circ\pi_{\frac{1}{n}, \frac{2}{n}, \dots, 1}$ \ is also measurable.

Next we show the measurability of \ $\Phi$.
Since the natural projection \ $\DD(\RR_+, \RR^d)\ni f\mapsto f(1)=\pi_1(f)$ \ is measurable,
 it is enough to show that the mapping
 \[
   \DD(\RR_+, \RR^d)\ni f\mapsto \widetilde\Phi(f):=\int_0^1 K( t, f(t) ) \, \dd t
 \]
 is measurable.
Namely, we show that \ $\widetilde\Phi$ \ is continuous.
We have to check that \ $\widetilde\Phi(f_n) \to \widetilde\Phi(f)$ \ in \ $\RR^{p}$ \ as
 \ $n \to \infty$ \ whenever \ $f_n \to f$ \ in \ $\DD(\RR_+, \RR^d)$ \ as
 \ $n \to \infty$, \ where \ $f,f_n\in D(\RR_+,\RR^d)$, $n\in\NN$.
\ Due to Ethier and Kurtz \cite[Proposition 3.5.3]{EthKur}, for all \ $T>0$ \
 there exists a sequence \ $\lambda_n : \RR_+ \to \RR_+$, \ $n \in \NN$, \ of strictly increasing
 continuous functions with \ $\lambda_n(0) = 0$ \ and
 \ $\lim_{t \to \infty} \lambda_n(t) = \infty$ \ such that
 \begin{align}\label{01seged15}
   \lim_{n \to \infty} \sup_{t \in [0, T]} | \lambda_n(t) - t | = 0 , \qquad
   \lim_{n \to \infty} \sup_{t \in [0, T]} \| f_n(t) - f(\lambda_n(t)) \| = 0 .
 \end{align}
We check that \ $\lim_{n\to\infty}f_n(t)=f(t)$, \ if \ $t\in\RR_+$ \ is a continuity point of \ $f$.
\ This readily follows by
 \[
    \Vert f_n(t)-f(t)\Vert
     \leq \Vert f_n(t)-f(\lambda_n(t))\Vert
           + \Vert f(\lambda_n(t))-f(t)\Vert,
      \qquad n\in\NN, \;\; t\in\RR_+.
 \]
Using that \ $f$ \ has at most countably many discontinuities (see, e.g., Jacod and Shiryaev
 \cite[page 326]{JSh}), we have \ $\lim_{n\to\infty}f_n(t)=f(t)$ \ for all \ $t\in\RR_+$ \ except
 a countable set having Lebesgue measure zero.
In what follows we check that
 \[
   \sup_{n \in \NN} \sup_{t \in [0, 1]} \| K( t, f_n(t) ) \| <\infty.
 \]
Since \ $K$ \ is continuous and hence it is bounded on a compact set, it is enough to verify that
 \[
    \sup_{n \in \NN} \sup_{t \in [0,1]} \| f_n(t) \|
       <\infty.
 \]
This follows by Jacod and Shiryaev \cite[Chapter VI, Lemma 1.14 (b)]{JSh},
 since \ $f_n\to f$ \ in \ $D(\RR_+,\RR^d)$ \ yields that \ $\{f_n : n\in\NN\}$ \ is a relatively
 compact set (with respect to the Skorokhod topology).
Then Lebesgue dominated convergence theorem yields the continuity of \ $\widetilde\Phi$.

In order to show \ $C_{\Phi,(\Phi_n)_{n \in \NN}} = \CC(\RR_+, \RR^d)$ \ we have to
 check that \ $\Phi_n(f_n) \to \Phi(f)$ \ whenever
 \ $f_n \lu f$ \ with \ $f \in \CC(\RR_+, \RR^d)$ \ and
 \ $f_n \in \DD(\RR_+, \RR^d)$, \ $n \in \NN$.
\ We have
 \begin{align*}
  \| \Phi_n(f_n) - \Phi(f) \|
  &\leq \| f_n(1) - f(1) \| +
      \frac{1}{n}
        \sum_{k = 1}^n
         \left\| K \left( \frac{k}{n}, f_n \left( \frac{k}{n} \right) \right)
                 - K \left( \frac{k}{n}, f \left( \frac{k}{n} \right) \right)
         \right\| \\
  &\quad + \sum_{k = 1}^n
            \int_{(k-1)/n}^{k/n}
             \left\| K \left( \frac{k}{n}, f \left( \frac{k}{n} \right) \right)
                     - K( t, f(t) ) \right\| \, \dd t\\
  &=: \| f_n(1) - f(1) \| + A_n^{(1)} + A_n^{(2)}.
 \end{align*}
Since \ $f_n \lu f$, \ we get
 \[
    \| f_n(1) - f(1) \|
      \leq \sup_{t \in [0,1]} \| f_n(t) - f(t)\|\to 0
      \qquad \text{as \ $n\to\infty$.}
 \]
Let us also observe that
 \[
   \sup_{n \in \NN} \sup_{t \in [0, 1]} \| f_n(t) \|
   \leq \sup_{n \in \NN} \sup_{t \in [0, 1]}
         \| f_n(t) - f(t) \| + \sup_{t \in [0, 1]} \| f(t) \|
   =: c < \infty ,
 \]
 hence
 \[
   A_n^{(1)} \leq C_{c} \sup_{t \in [0, 1]} \| f_n(t) - f(t) \| \to 0
 \]
 as \ $n \to \infty$.
\ Moreover,
 \[
   A_n^{(2)}
   \leq C_{c}
        \sum_{k = 1}^n
         \int_{(k-1)/n}^{k/n}
          \left( \left| \frac{k}{n} - t \right|
                 + \left\| f \left( \frac{k}{n} \right) - f(t) \right\|
          \right) \dd t
   \leq C_{c} ( n^{-1} + \omega_1(f, n^{-1}) ) ,
 \]
 where
 \[
   \omega_1(f, \vare)
   := \sup_{t, \, s \in [0,1], \, |t-s|<\vare}
       \|\ f(t) - f(s) \| , \qquad \vare > 0,
 \]
 denotes the modulus of continuity of \ $f$ \ on \ $[0, 1]$.
\ Since \ $f$ \ is continuous, \ $\omega_1(f, n^{-1}) \to 0$ \ as
 \ $n \to \infty$ \ (see, e.g., Jacod and Shiryaev
 \cite[Chapter VI, 1.6]{JSh}), and we obtain \ $A_n^{(2)}\to0$ \ as \ $n\to\infty$.
Then \ $C_{\Phi,(\Phi_n)_{n \in \NN}} = \CC(\RR_+, \RR^d)$.

Finally, \ $\CC(\RR_+, \RR^d) \in \cB(\DD(\RR_+, \RR^d))$ \ holds since
 \ $\DD(\RR_+, \RR^d) \setminus \CC(\RR_+, \RR^d)$ \ is open.
Indeed, if \ $f \in \DD(\RR_+, \RR^d) \setminus \CC(\RR_+, \RR^d)$ \ then
 there exists \ $t \in \RR_+$ \ such that
 \ $\vare := \| f(t) - \lim_{s \uparrow t} f(s) \| > 0$, \ and then the open
 ball in \ $\DD(\RR_+, \RR^d)$ \ with centre \ $f$ \ and radius \ $\vare / 2$
 \ does not contain any continuous function.
We note that for \ $\CC(\RR_+, \RR^d) \in \cB(\DD(\RR_+, \RR^d))$ \ one can also simply
 refer to Ethier and Kurtz \cite[Problem 3.11.25]{EthKur}.
\proofend

The next lemma is a key tool for proving joint convergence in distribution of
 \ $n^{-5}\det(\bA_n)$ \ and \ $n^{-4}\tbA_n \bd_n$.
\ We collected all the 'building blocks' that will appear in the asymptotic expansions
 of \ $\det(\bA_n)$, \ $\tbA_n$ \ and \ $\bd_n$.

\begin{Lem}\label{Lem01_UV3}
Let \ $(X_k)_{k \geq -1}$ \ be a nonprimitive INAR(2) process with autoregressive parameter
  \ $(0,1)$.
\ Suppose that \ $X_0 = X_{-1} = 0$, \ $\EE(\vare_1^2) < \infty$ \ and
 \ $\mu_\vare > 0$.
\ Then
  \begin{align}\label{01seged_main}
   (S_n^{(i)})_{i=1}^9 \distr (S^{(i)})_{i=1}^9 \qquad \text{as \ $n \to \infty$,}
 \end{align}
 where, for all \ $n\in\NN$,
 \begin{alignat*}{2}
  S_n^{(1)} &:= \frac{1}{n^{1/2}} ( U_n - \EE(U_n) ), \qquad&
  S_n^{(2)} &:= \frac{1}{n^{1/2}} ( V_n - \EE(V_n) ), \\
  S_n^{(3)} &:= \frac{1}{n^{5/2}} \sum_{k=1}^n k ( U_k - \EE(U_k) ), \qquad&
  S_n^{(4)} &:= \frac{1}{n^{5/2}} \sum_{k=1}^n k ( V_k - \EE(V_k) ), \\
  S_n^{(5)} &:= \frac{1}{n^2} \sum_{k=1}^n ( U_k - \EE(U_k) )^2, \qquad&
  S_n^{(6)} &:= \frac{1}{n^2} \sum_{k=1}^n ( V_k - \EE(V_k) )^2, \\
  S_n^{(7)} &:= \frac{1}{n^{3/2}} \sum_{k=1}^n ( U_k - \EE(U_k) ), \qquad&
  S_n^{(8)} &:= \frac{1}{n^{3/2}} \sum_{k=1}^n ( V_k - \EE(V_k) ),
 \end{alignat*}
 \[
   S_n^{(9)} := \frac{1}{n^2}
                \sum_{k=1}^n ( U_k - \EE(U_k) )( V_k - \EE(V_k) ),
 \]
 and
 \begin{alignat*}{3}
  S^{(1)} &:= \sigma_\vare \cW_1^{(1)} , \qquad&
  S^{(2)} &:= \sigma_\vare \cW_1^{(2)} , \qquad&
  S^{(3)} &:= \sigma_\vare \int_0^1 t \cW_t^{(1)} \, \dd t , \\
  S^{(4)} &:= \sigma_\vare \int_0^1 t \cW_t^{(2)} \, \dd t , \qquad&
  S^{(5)} &:= \sigma_\vare^2 \int_0^1 (\cW_t^{(1)})^2 \, \dd t , \qquad&
  S^{(6)} &:= \sigma_\vare^2 \int_0^1 (\cW_t^{(2)})^2 \, \dd t , \\
  S^{(7)} &:= \sigma_\vare \int_0^1 \cW_t^{(1)} \, \dd t , \qquad&
  S^{(8)} &:= \sigma_\vare \int_0^1 \cW_t^{(2)} \, \dd t ,  \qquad&
  S^{(9)} &:= \sigma_\vare^2 \int_0^1 \cW_t^{(1)} \cW_t^{(2)} \, \dd t ,
 \end{alignat*}
 where \ $(\cW_t^{(1)})_{t \in \RR_+}$ \ and \ $(\cW_t^{(2)})_{t \in \RR_+}$ \ are
 independent standard Wiener processes.
Especially, for all \ $\delta > 0$,
 \begin{align}\label{01seged10}
   \frac{1}{n^\delta} (S_n^{(i)})_{i=1}^9 \stoch {\bf 0}\in\RR^9 \qquad
   \text{as \ $n \to \infty$.}
 \end{align}
\end{Lem}

\noindent{\bf Proof.}
The proof is based on Lemma \ref{Conv2Funct} and Lemma \ref{Lem01_UV2}.
Let us introduce the function \ $K : [0,1] \times \RR^2 \to \RR^7$ \ defined by
 \[
   K( t, (x_1, x_2) ) := (tx_1, tx_2, x_1^2, x_2^2, x_1, x_2, x_1x_2) ,
   \qquad t \in [0,1] , \quad (x_1, x_2) \in \RR^2 .
 \]
Then \eqref{01seged5} holds, since for all \ $s,t\in[0,1]$, \ $R>0$ \ and
 \ $x=(x_1,x_2)\in\RR^2$, \ $y=(y_1,y_2)\in\RR^2$ \ with \ $\Vert x\Vert\leq R$
 \ and \ $\Vert y\Vert\leq R$, \ we get
 \begin{align*}
   &\Vert K(t,x_1,x_2) - K(s,y_1,y_2)\Vert\\
   &\quad = \Vert (tx_1-sy_1,tx_2-sy_2,x_1^2-y_1^2,x_2^2-y_2^2,x_1-y_1,x_2-y_2,x_1x_2-y_1y_2) \Vert \\
   &\quad  \leq \Big( 2t^2(x_1-y_1)^2 + 2y_1^2(t-s)^2 + 2t^2(x_2-y_2)^2
                + 2y_2^2(t-s)^2 + (x_1-y_1)^2(x_1+y_1)^2 \\
   &\phantom{\quad  \leq \Big( }
                + (x_2-y_2)^2(x_2+y_2)^2  + (x_1-y_1)^2 + (x_2-y_2)^2 + 2x_2^2(x_1-y_1)^2
                + 2y_1^2(x_2-y_2)^2\Big)^{1/2}\\
   &\quad \leq \max(\sqrt{2},2R)\big(4(x_1-y_1)^2 + 4(x_2-y_2)^2 + 2(t-s)^2\big)^{1/2}\\
   &\quad \leq 2\max(\sqrt{2},2R)(\vert t-s\vert + \Vert x-y\Vert),
 \end{align*}
 where the last step follows by Minkowski's inequality.
Further,
 \begin{align*}
  \Phi_n( n^{-1/2} ( \cU^n - \EE(\cU^n) ), n^{-1/2} ( \cV^n - \EE(\cV^n) ) )
  &= (S_n^{(i)})_{i=1}^9 , \qquad n \in \NN , \\
  \Phi( \sigma_\vare \cW^{(1)}, \sigma_\vare \cW^{(2)} )
  &= (S^{(i)})_{i=1}^9 ,
 \end{align*}
 where \ $(\Phi_n)_{n \in \NN}$ \ and \ $\Phi$ \ are defined in Lemma
 \ref{Lem01_UV2}.
By Lemma \ref{Lem01_UV2}, \ $C_{\Phi,(\Phi_n)_{n \in \NN}} = \CC(\RR_+, \RR^2)\in\cB(\DD(\RR_+, \RR^2))$
 \ and using that a standard Wiener process has continuous trajectories with probability one, we have
 \ $\PP((\cW^{(1)},\cW^{(2)})\in C_{\Phi,(\Phi_n)_{n \in \NN}} )=1$. \
Since \ $\CC(\RR_+, \RR^2)$ \ is a measurable subset of \ $\DD(\RR_+, \RR^2)$ \ (see  Lemma \ref{Lem01_UV2}),
 Lemma \ref{Conv2Funct} and Lemma \ref{Lem01_UV} imply \eqref{01seged_main}.

Finally, Slutsky's lemma yields \eqref{01seged10}.
\proofend

The next lemma describes the asymptotic behavior of the sequence \ $(\det({\bf A}_n))_{n\in\NN}$.

\begin{Lem}\label{Lem01}
Let \ $(X_k)_{k \geq -1}$ \ be a nonprimitive \INARtwo\ process with autoregressive parameter
 \ $(0,1)$.
\ Suppose that \ $X_0 = X_{-1} = 0$, \ $\EE(\vare_1^2) < \infty$ \ and \ $\mu_\vare > 0$.
\ Then
 \begin{equation}\label{01det}
  n^{-5} \det(\bA_n)
  \distr \frac{\mu_\vare^2 \sigma_\vare^2}{12} \int_0^1 (\cW_t)^2 \, \dd t
  \qquad \text{as \ $n\to\infty$,}
 \end{equation}
 where \ $(\cW_t)_{t \in \RR_+}$ \ is a standard Wiener process.
\end{Lem}

\noindent
\textbf{Proof.}
In order to investigate the asymptotic behavior of the sequence
 \ $(\det(\bA_n))_{n \in \NN}$, \ we derive (asymptotic) expansions for the entries
 of the matrices \ $\bA_n$, $n\in\NN$.
\ First we separate the expectations (`leading terms') of entries of \ $\bA_n$.
\ Namely, we get
 \begin{align*}
  \sum_{k=1}^{2n} X_{k-1}^2
  &= \sum_{k=1}^{n-1} U_k^2 + \sum_{k=1}^n V_k^2
   = \sum_{k=1}^{n-1}
      \Big[ k \mu_\vare + \big( U_k - \EE(U_k) \big) \Big]^2
     + \sum_{k=1}^n
        \Big[ k \mu_\vare + \big( V_k - \EE(V_k) \big) \Big]^2 \\
  &= \mu_\vare^2 \left[ \sum_{k=1}^{n-1} k^2 + \sum_{k=1}^n k^2 \right]
     + 2 \mu_\vare \left[ \sum_{k=1}^{n-1} k \big( U_k - \EE(U_k) \big)
                         + \sum_{k=1}^n k \big( V_k - \EE(V_k) \big) \right] \\
&\phantom{\quad}
     +  \sum_{k=1}^{n-1} \big( U_k - \EE(U_k) \big)^2
     +  \sum_{k=1}^n \big( V_k - \EE(V_k) \big)^2  \\
   &=  F_{2n,1} (2n)^3 + F_{2n,2} (2n)^{5/2} + F_{2n,3} (2n)^2 ,
 \end{align*}
 where, by \eqref{01seged_main},
 \begin{align*}
  F_{2n,1} &:= \frac{\mu_\vare^2}{(2n)^3}
              \left[ \sum_{k=1}^{n-1} k^2 + \sum_{k=1}^n k^2 \right]
           = \frac{(2n^2+1)\mu_\vare^2}{24n^2}
           \to\frac{\mu_\vare^2}{12} =: F_1 , \\[2mm]
  F_{2n,2}&:= \frac{2\mu_\vare}{(2n)^{5/2}}
            \left[ \sum_{k=1}^{n-1} k \big( U_k - \EE(U_k) \big)
                  + \sum_{k=1}^n k \big( V_k - \EE(V_k) \big) \right] \\
          &\distr
           \frac{\mu_\vare\sigma_\vare}{2^{3/2}}
           \int_0^1 s \big(\cW^{(1)}_s + \cW^{(2)}_s \big) \,\dd s =: F_2 , \\[2mm]
  F_{2n,3}&:= \frac{1}{(2n)^2}
               \left[\sum_{k=1}^{n-1} \big( U_k - \EE(U_k) \big)^2
               +  \sum_{k=1}^n \big( V_k - \EE(V_k) \big)^2
               \right] \\
          &\distr
            \frac{\sigma_\vare^2}{4}
            \int_0^1 \Big[ \big( \cW^{(1)}_s \big)^2
                          + \big( \cW^{(2)}_s \big)^2 \Big] \dd s = :F_3 ,
 \end{align*}
 where \ $(\cW^{(1)}_t)_{t\in\RR_+}$ \ and \ $(\cW^{(2)}_t)_{t\in\RR_+}$ \ are independent
 standard Wiener processes.
In a similar way
 \begin{align*}
  \sum_{k=1}^{2n} X_{k-2}^2
  &= \sum_{k=1}^{n-1} U_k^2 + \sum_{k=1}^{n-1} V_k^2
   = \sum_{k=1}^{n-1}
      \Big[ k \mu_\vare + \big( U_k - \EE(U_k) \big) \Big]^2
     + \sum_{k=1}^{n-1}
        \Big[ k \mu_\vare + \big( V_k - \EE(V_k) \big) \Big]^2 \\
  &= 2\mu_\vare^2 \sum_{k=1}^{n-1} k^2
     + 2 \mu_\vare \sum_{k=1}^{n-1} k \big( U_k - \EE(U_k) + V_k - \EE(V_k) \big) \\
&\phantom{\quad}
     + \sum_{k=1}^{n-1} \big( U_k - \EE(U_k) \big)^2
     + \sum_{k=1}^{n-1} \big( V_k - \EE(V_k) \big)^2 \\
   &=  H_{2n,1} (2n)^3 + H_{2n,2} (2n)^{5/2} + H_{2n,3} (2n)^2 ,
 \end{align*}
 where, by \eqref{01seged_main}, Slutsky's lemma and continuity theorem,
 \begin{align*}
  H_{2n,1} &:= \frac{2\mu_\vare^2}{(2n)^3} \sum_{k=1}^{n-1} k^2
           = \frac{(n-1)(2n-1)\mu_\vare^2}{24n^2}
           \to \frac{\mu_\vare^2}{12} =: H_1=F_1 , \\[2mm]
  H_{2n,2}&:= \frac{2\mu_\vare}{(2n)^{5/2}}
             \sum_{k=1}^{n-1} k \big( U_k - \EE(U_k) + V_k - \EE(V_k) \big) \\
          &\distr
           \frac{\mu_\vare\sigma_\vare}{2^{3/2}}
           \int_0^1 s \big(\cW^{(1)}_s + \cW^{(2)}_s \big) \,\dd s =: H_2=F_2 , \\[2mm]
  H_{2n,3}&:= \frac{1}{(2n)^2}
             \left[\sum_{k=1}^{n-1} \big( U_k - \EE(U_k) \big)^2
               +  \sum_{k=1}^{n-1} \big( V_k - \EE(V_k) \big)^2\right] \\
          &\distr
            \frac{\sigma_\vare^2}{4}
            \int_0^1 \Big[ \big( \cW^{(1)}_s \big)^2
                          + \big( \cW^{(2)}_s \big)^2 \Big] \dd s = :H_3=F_3.
 \end{align*}
Further,
 \begin{align*}
  \sum_{k=1}^{2n} X_{k-1} X_{k-2}
  &= \sum_{k=1}^{n-1} U_k V_k + \sum_{k=1}^{n-1} U_k V_{k+1} \\
  &= \sum_{k=1}^{n-1}
      \Big[ k \mu_\vare + \big( U_k - \EE(U_k) \big) \Big]
      \Big[ (2k+1) \mu_\vare + \big( V_k - \EE(V_k) \big)
                            + \big( V_{k+1} - \EE(V_{k+1}) \big) \Big] \\
  &= \mu_\vare^2 \sum_{k=1}^{n-1} k(2k+1)
     + \sum_{k=1}^{n-1}
        \big( U_k - \EE(U_k) \big)
        \big( V_k - \EE(V_k) + V_{k+1} - \EE(V_{k+1}) \big) \\
&\phantom{\quad}
     + \mu_\vare
       \sum_{k=1}^{n-1} \Big[ k \big( V_k - \EE(V_k) + V_{k+1} - \EE(V_{k+1}) \big)
                            + (2k+1) \big( U_k - \EE(U_k) \big) \Big] \\
   &=  G_{2n,1} (2n)^3 + G_{2n,2} (2n)^{5/2} + G_{2n,3} (2n)^2 ,
 \end{align*}
 where, by \eqref{01seged_main}, \eqref{01seged10}, Slutsky's lemma and continuity theorem,
 \begin{align*}
  G_{2n,1} &:= \frac{\mu_\vare^2}{(2n)^3} \sum_{k=1}^{n-1} k(2k+1)
           = \frac{(4n+1)(n-1)\mu_\vare^2}{48n^2}
           \to \frac{\mu_\vare^2}{12} =: G_1=F_1 , \\[2mm]
  G_{2n,2}&:= \frac{\mu_\vare}{(2n)^{5/2}}
             \sum_{k=1}^{n-1}
              \Big[ k \big( V_k - \EE(V_k) + V_{k+1} - \EE(V_{k+1}) \big)
                    + (2k+1) \big( U_k - \EE(U_k) \big) \Big]\\
          &\distr
           \frac{\mu_\vare\sigma_\vare}{2^{3/2}}
           \int_0^1 s \big(\cW^{(1)}_s + \cW^{(2)}_s \big) \,\dd s =: G_2=F_2 , \\[2mm]
  G_{2n,3}&:= \frac{1}{(2n)^2}
             \sum_{k=1}^{n-1}
              \big( U_k - \EE(U_k) \big)
              \big( V_k - \EE(V_k) + V_{k+1} - \EE(V_{k+1}) \big)\\
          &\distr
            \frac{\sigma_\vare^2}{2}
            \int_0^1 \cW^{(1)}_s \cW^{(2)}_s \, \dd s = :G_3 .
 \end{align*}
For the derivation of the convergence in distribution \ $G_{2n,3}\distr G_3$ \ as \ $n\to\infty$,
 \ we give a bit more explanation.
Slutsky's lemma, \eqref{01seged_main} and \eqref{01seged10} follow the desired convergence if we check that
 \[
   \frac{1}{n^2}\sum_{k=1}^{n-1} (U_k - \EE(U_k))(\vare_{2k+1} - \mu_\vare)
     \stoch 0
   \qquad \text{as \ $n\to\infty$.}
 \]
In fact, we prove that the above convergence holds in $L_1$-sense.
Namely, by Cauchy-Schwarz's inequality, we get
 \begin{align*}
  &\EE\left\vert  \frac{1}{n^2}\sum_{k=1}^{n-1} (U_k - \EE(U_k))(\vare_{2k+1} - \mu_\vare) \right\vert
      \leq \frac{1}{n^2}\sum_{k=1}^{n-1}\EE\vert (U_k - \EE(U_k))(\vare_{2k+1} - \mu_\vare)\vert\\
  &\qquad \leq \frac{1}{n^2}\sum_{k=1}^{n-1} \sqrt{\EE(U_k - \EE(U_k))^2 \EE(\vare_{2k+1} - \mu_\vare)^2}
       =\frac{1}{n^2}\sum_{k=1}^{n-1}\sqrt{k\sigma_\vare^2\cdot\sigma_\vare^2} \\
  &\qquad  =\frac{\sigma_\vare^2}{n^2}\sum_{k=1}^{n-1}\sqrt{k}
       \leq \frac{\sigma_\vare^2}{n^2}(n-1)^{3/2} \to 0
       \qquad \text{as \ $n\to\infty$.}
 \end{align*}

Similar expansions can be derived for \ $\sum_{k=1}^{2n+1} X_{k-1}^2$,
 \ $\sum_{k=1}^{2n+1} X_{k-2}^2$ \ and \ $\sum_{k=1}^{2n+1} X_{k-1} X_{k-2}$.
\ Namely,
  \begin{align*}
     &\sum_{k=1}^{2n+1} X_{k-1}^2
       =  F_{2n+1,1} (2n+1)^3 + F_{2n+1,2} (2n+1)^{5/2} + F_{2n+1,3} (2n+1)^2,
       \qquad n\in\ZZ_+, \\
     &  \sum_{k=1}^{2n+1} X_{k-2}^2
       =  H_{2n+1,1} (2n+1)^3 + H_{2n+1,2} (2n+1)^{5/2} + H_{2n+1,3} (2n+1)^2,
      \qquad n\in\ZZ_+, \\
     & \sum_{k=1}^{2n+1} X_{k-1}X_{k-2}
       =  G_{2n+1,1} (2n+1)^3 + G_{2n+1,2} (2n+1)^{5/2} + G_{2n+1,3} (2n+1)^2,
      \qquad n\in\ZZ_+,
  \end{align*}
 where
 \begin{align*}
   & F_{2n+1,1} :=\left(\frac{2n}{2n+1}\right)^3
                  \left(F_{2n,1} + \frac{n^2\mu_\vare^2}{(2n)^3}\right)
                =\left(\frac{2n}{2n+1}\right)^3 \frac{(n+1)(2n+1)\mu_\vare^2}{24n^2}\to F_1,\\[1mm]
   & F_{2n+1,2} := \left(\frac{2n}{2n+1}\right)^{5/2}
                  \left(F_{2n,2} + \frac{2\mu_\vare n(U_n - \EE(U_n))}{(2n)^{5/2}}\right)\\
   &\phantom{F_{2n+1,2} \;\;}
               =  \left(\frac{2n}{2n+1}\right)^{5/2}\frac{2\mu_\vare}{(2n)^{5/2}}\left[ \sum_{k=1}^n k \big( U_k - \EE(U_k) \big)
                   + \sum_{k=1}^n k \big( V_k - \EE(V_k) \big) \right] \distr F_2,\\[1mm]
   & F_{2n+1,3} :=  \left(\frac{2n}{2n+1}\right)^2\left(F_{2n,3} + \frac{(U_n - \EE(U_n))^2}{(2n)^2}\right)\\
   &\phantom{F_{2n+1,3} \;\;} =
                 \left(\frac{2n}{2n+1}\right)^2\frac{1}{(2n)^2}
                   \left[\sum_{k=1}^n \big( U_k - \EE(U_k) \big)^2
                    + \sum_{k=1}^n \big( V_k - \EE(V_k) \big)^2
                    \right]\distr F_3,
 \end{align*}
 and
 \begin{align*}
   & H_{2n+1,1} := \left(\frac{2n}{2n+1}\right)^3 \left( H_{2n,1} + \frac{n^2\mu_\vare^2}{(2n)^3}\right)
                 = \left(\frac{2n}{2n+1}\right)^3  \frac{(2n^2+1)\mu_\vare^2}{24n^2}\to H_1,\\
   & H_{2n+1,2} := \left(\frac{2n}{2n+1}\right)^{5/2} \left(H_{2n,2} + \frac{2\mu_\vare n(V_n - \EE(V_n))}{(2n)^{5/2}}\right)\\
   &\phantom{H_{2n+1,2} \;\;}
               = \left(\frac{2n}{2n+1}\right)^{5/2} \frac{2\mu_\vare}{(2n)^{5/2}}\left[ \sum_{k=1}^{n-1} k \big( U_k - \EE(U_k) \big)
                   + \sum_{k=1}^n k \big( V_k - \EE(V_k) \big) \right] \distr H_2,\\
   & H_{2n+1,3} := \left(\frac{2n}{2n+1}\right)^2 \left(H_{2n,3} + \frac{(V_n - \EE(V_n))^2}{(2n)^2}\right)\\
   &\phantom{H_{2n+1,3} \;\;}
                 = \left(\frac{2n}{2n+1}\right)^2 \frac{1}{(2n)^2}
                   \left[\sum_{k=1}^{n-1} \big( U_k - \EE(U_k) \big)^2
                    + \sum_{k=1}^n \big( V_k - \EE(V_k) \big)^2
                    \right]\distr H_3,
 \end{align*}
 and
 \begin{align*}
   & G_{2n+1,1} := \left(\frac{2n}{2n+1}\right)^3 \left(G_{2n,1} + \frac{n^2\mu_\vare^2}{(2n)^3}\right)
                 = \left(\frac{2n}{2n+1}\right)^3
                     \frac{(4n-1)(n+1)\mu_\vare^2}{48n^2}\to G_1,\\
   & G_{2n+1,2} :=\left(\frac{2n}{2n+1}\right)^{5/2}\left(G_{2n,2} + \frac{n\mu_\vare ( U_n - \EE(U_n) + V_n - \EE(V_n) )}{(2n)^{5/2}}\right)
                                \distr G_2,\\
   & G_{2n+1,3} := \left(\frac{2n}{2n+1}\right)^2\left( G_{2n,3} + \frac{(U_n - \EE(U_n))(V_n - \EE(V_n))}{(2n)^2}\right) \distr G_3.
 \end{align*}
Hence we have an (asymptotic) expansion for \ $\det(\bA_n)$, \ namely,
 \begin{align} \label{A}
  \begin{split}
   \det(\bA_n)
   &= ( F_{n,1} n^3 + F_{n,2} n^{5/2} + F_{n,3} n^2 )
      ( H_{n,1} n^3 + H_{n,2} n^{5/2} + H_{n,3} n^2 ) \\
   &\phantom{\quad}
      - ( G_{n,1} n^3 + G_{n,2} n^{5/2} + G_{n,3} n^2 )^2 \\
   &= (F_{n,1} H_{n,1} - G_{n,1}^2) n^6
      + (F_{n,1} H_{n,2} - 2 G_{n,1} G_{n,2} + F_{n,2} H_{n,1}) n^{11/2} \\
   &\phantom{\quad}
      + (F_{n,1} H_{n,3} + F_{n,2} H_{n,2} + F_{n,3} H_{n,1} - 2 G_{n,1} G_{n,3}
        - G_{n,2}^2) n^5 \\
   &\phantom{\quad}
      + (F_{n,2} H_{n,3} - 2 G_{n,2} G_{n,3} + F_{n,3} H_{n,2}) n^{9/2}
      + (F_{n,3} H_{n,3} - G_{n,3}^2) n^4 ,
      \qquad n\in\NN.
  \end{split}
 \end{align}
By Lemma \ref{Lem01_UV3}, Slutsky's lemma and the continuous mapping theorem, we have
 \ $F_{n,i},G_{n,i}, H_{n,i}$, $i=1,2,3$, \ converge jointly in distribution,
 and hence the coefficients of the expansion \eqref{A} also converge jointly in distribution.
Futher, we show that the first two leading terms have no influence by which we mean that
 \begin{gather}
  n (F_{n,1}H_{n,1} - G_{n,1}^2) \to 0 , \label{FHG2} \\
  n^{1/2} (F_{n,1}H_{n,2} - 2G_{n,1}G_{n,2} + F_{n,2}H_{n,1}) \stoch 0 , \label{FHGGFH}
 \end{gather}
 as \ $n\to\infty$.
\ Indeed,
 \begin{align*}
  2n (F_{2n,1}H_{2n,1} - G_{2n,1}^2)
    &=\frac{2\mu_\vare^4}{(48)^2 n^3}
       \Big[ 4 (2n^2 + 1)(n-1)(2n-1) - (4n^2 - 3n - 1)^2 \Big] \\
    &=\frac{2\mu_\vare^4(n-1)(15n-3)}{(48)^2 n^3} \to 0
    \qquad \text{as \ $n\to\infty$.}
 \end{align*}
We note that \ $(2n+1)(F_{2n+1,1}H_{2n+1,1} - G_{2n+1,1}^2)\to 0$ \ as \ $n\to\infty$
 \ can be proved similarly.
Indeed,
 \begin{align*}
   &(2n+1)(F_{2n+1,1}H_{2n+1,1} - G_{2n+1,1}^2)\\
    & =\frac{(2n)^6}{(2n+1)^5}\left( \left(F_{2n,1} + \frac{\mu_\vare^2}{8n}\right)
                   \left(H_{2n,1} + \frac{\mu_\vare^2}{8n}\right)
                   - \left(G_{2n,1} +\frac{\mu_\vare^2}{8n}\right)^2
            \right)\\
    & = \frac{(2n)^6}{(2n+1)^5} \left(F_{2n,1}H_{2n,1} - G_{2n,1}^2
                   + \frac{\mu_\vare^2}{8n}(F_{2n,1} + H_{2n,1})
                   - \frac{\mu_\vare^2}{4n}G_{2n,1}
           \right) \\
    & = \frac{(2n)^6}{(2n+1)^5}\left( \left(F_{2n,1}H_{2n,1} - G_{2n,1}^2 \right)
       + \mu_\vare^4\frac{4n^2-3n+2}{192 n^3}
       - \mu_\vare^4\frac{4n^2-3n-1}{192 n^3} \right) \\
    & =\left(\frac{2n}{2n+1}\right)^5 2n\left(F_{2n,1}H_{2n,1} - G_{2n,1}^2 \right)
        + \frac{\mu_\vare^4 n^3}{8(2n+1)^5}
     \to 0
     \qquad \text{as \ $n\to\infty$.}
 \end{align*}
Now we turn to check \eqref{FHGGFH}.
We get
 \begin{align*}
  &(2n)^{1/2} (F_{2n,1}H_{2n,2} - 2G_{2n,1}G_{2n,2} + F_{2n,2}H_{2n,1}) \\
  &= \frac{\mu_\vare^3}{96 n^4}
     \left[ 2(2 n^2 + 1) \sum_{k=1}^{n-1} k (U_k - \EE(U_k) + V_k - \EE(V_k))
     \right. \\
  &\phantom{\quad \frac{\mu_\vare^3}{96 n^4} \bigg[}
            + 2(n-1)(2n-1)
              \left( \sum_{k=1}^{n-1} k (U_k - \EE(U_k))
                      + \sum_{k=1}^n k (V_k - \EE(V_k)) \right) \\
  &\phantom{\quad \frac{\mu_\vare^3}{96 n^4} \bigg[}
     \left. - (4n+1)(n-1)
       \sum_{k=1}^{n-1}
        \Big( k\big(V_k - \EE(V_k) + V_{k+1} - \EE(V_{k+1})\big)
               + (2k+1) \big(U_k - \EE(U_k)\big) \Big) \right] \\
  &= \frac{\mu_\vare^3}{16 n^4}
     \sum_{k=1}^{n-1} k (U_k - \EE(U_k) + V_k - \EE(V_k))
     - \frac{(4n+1)(n-1)\mu_\vare^3}{96 n^4}
       \sum_{k=1}^{n-1} (U_k - \EE(U_k)) \\
  &\phantom{\quad}
     + \frac{(4n+1)(n-1)\mu_\vare^3}{96 n^4} (V_1 - \EE(V_1))
     + \frac{(1-n)\mu_\vare^3}{32 n^3} (V_n - \EE(V_n)) \\
  &\phantom{\quad}
    + \frac{(4n+1)(n-1)\mu_\vare^3}{96 n^4}
        \sum_{k=2}^n (V_k-\EE(V_k))
   \stoch0,
 \end{align*}
 where we used \eqref{01seged10} and that
 \ $n^{-\delta}(V_1 - \EE(V_1))\as 0$ \ as \ $n\to\infty$ \ for all \ $\delta>0$.
\ Similarly one can prove that
 \[
  (2n+1)^{1/2} (F_{2n+1,1}H_{2n+1,2} - 2G_{2n+1,1}G_{2n+1,2} + F_{2n+1,2}H_{2n+1,1})\stoch 0
  \qquad\text{as \ $n\to\infty$.}
 \]
Indeed,
 \begin{align*}
 (2n+1)^{1/2} &(F_{2n+1,1}H_{2n+1,2} - 2G_{2n+1,1}G_{2n+1,2} + F_{2n+1,2}H_{2n+1,1}) \\
   &  = \frac{(2n)^{11/2}}{(2n+1)^5}
      \Bigg[  \left( F_{2n,1}+\frac{\mu_\vare^2}{8n}\right)
              \left( H_{2n,2}+ \mu_\vare\frac{V_n - \EE(V_n)}{(2n)^{3/2}}\right)  \\
  & \phantom{ = (2n+1)^{1/2} \Bigg[}
            - 2 \left( G_{2n,1}+\frac{\mu_\vare^2}{8n}\right)
                \left( G_{2n,2}+ \mu_\vare\frac{U_n - \EE(U_n) + V_n - \EE(V_n)}{2^{5/2} n^{3/2}}\right) \\
  & \phantom{ = (2n+1)^{1/2} \Bigg[}
           +  \left( F_{2n,2}+ \mu_\vare\frac{U_n - \EE(U_n)}{(2n)^{3/2}} \right)
              \left( H_{2n,1}+ \frac{\mu_\vare^2}{8n} \right)
      \Bigg]\\
   & = \frac{(2n)^{11/2}}{(2n+1)^5}
         \big[(F_{2n,1}H_{2n,2} - 2G_{2n,1}G_{2n,2} + F_{2n,2}H_{2n,1})
               + R_n\big],
 \end{align*}
 where
 \begin{align*}
   R_n:=
     & F_{2n,1}\mu_\vare \frac{V_n-\EE(V_n)}{(2n)^{3/2}}
       + \frac{\mu_\vare^2}{8n}H_{2n,2}
       +  \frac{\mu_\vare^3}{8n}\cdot\frac{V_n-\EE(V_n)}{(2n)^{3/2}}
       - 2 G_{2n,1}\mu_\vare \frac{U_n-\EE(U_n) + V_n-\EE(V_n)}{2^{5/2} n^{3/2}} \\
    &  - \frac{\mu_\vare^2}{4n}G_{2n,2}
       - \frac{\mu_\vare^3}{4n}\cdot\frac{U_n-\EE(U_n) + V_n-\EE(V_n)}{2^{5/2} n^{3/2}}
       + F_{2n,2}\frac{\mu_\vare^2}{8n}
       + \mu_\vare \frac{U_n-\EE(U_n)}{(2n)^{3/2}} H_{2n,1}\\
     & + \frac{\mu_\vare^3}{8n}\cdot\frac{U_n-\EE(U_n)}{(2n)^{3/2}}.
 \end{align*}
Using that \ $F_{n,i}$, $G_{n,i}$, $H_{n,i}$, $i=1,2,3$, \ converge jointly in distribution,
 by \eqref{01seged10}, we get \ $(2n+1)^{1/2} R_n\stoch 0$ \ as \ $n\to\infty$, \ and hence
 Slutsky's lemma yields the desired convergence.

Using again the above mentioned joint convergence of \ $F_{n,i}$, $G_{n,i}$, $H_{n,i}$, $i=1,2,3$,
 \ \eqref{A}, \eqref{FHG2}, \eqref{FHGGFH} and Slutsky's lemma imply
 \begin{align}\label{01seged16}
  n^{-5} \det(\bA_n)
  \distr F_1 H_3 + F_2 H_2 + F_3 H_1 - 2 G_1 G_3 - G_2^2
  = \frac{\mu_\vare^2 \sigma_\vare^2}{24}
    \int_0^1 \big( \cW_t^{(1)} - \cW_t^{(2)} \big)^2 \dd t .
 \end{align}
It is easy to check that \ $2^{-1/2} \big( \cW_t^{(1)} - \cW_t^{(2)} \big)$,
 \ $t \in \RR_+$, \ is a standard Wiener process, hence the proof of
 \eqref{01det} is complete.
\proofend

\vskip0.5cm

\noindent
\textbf{Proof of Theorem \ref{01main}.} \
Using the (asymptotic) expansions derived in the proof of Lemma \ref{Lem01} for
 \ $\sum_{k=1}^n X_{k-1}^2$, \ $\sum_{k=1}^n X_{k-2}^2$ \ and
 \ $\sum_{k=1}^n X_{k-1} X_{k-2}$, \ we obtain an (asymptotic) expansion for the adjoint
 \ $\tbA_n$ \ of \ ${\bf A_n}$:
 \begin{align}\label{01seged2}
  \tbA_n = \tbA_{n,1} n^3 + \tbA_{n,2} n^{5/2} + \tbA_{n,3} n^2 ,
  \qquad n\in\NN,
 \end{align}
 where
 \begin{align*}
  \tbA_{n,1} &:= \begin{bmatrix}
               H_{n,1} & -G_{n,1} \\
               -G_{n,1} & F_{n,1}
              \end{bmatrix}
           \to \frac{\mu_\vare^2}{12}
               \begin{bmatrix}
                1 & -1 \\
                -1 & 1
               \end{bmatrix} =: \tbA^{(1)}, \\
  \tbA_{n,2} &:= \begin{bmatrix}
               H_{n,2} & -G_{n,2} \\
               -G_{n,2} & F_{n,2}
              \end{bmatrix}
           \distr \frac{\mu_\vare \sigma_\vare}{2^{3/2}}
                  \int_0^1 t \big( \cW_t^{(1)} + \cW_t^{(2)} \big) \dd t
                  \begin{bmatrix}
                   1 & -1 \\
                   -1 & 1
                  \end{bmatrix}  =: \tbA^{(2)}, \\
  \tbA_{n,3} &:= \begin{bmatrix}
               H_{n,3} & -G_{n,3} \\
               -G_{n,3} & F_{n,3}
              \end{bmatrix} \\
            &\distr \frac{\sigma_\vare^2}{4}
                  \begin{bmatrix}
                   \DS\int_0^1 \Big( \big( \cW^{(1)}_t \big)^2
                                  + \big( \cW^{(2)}_t \big)^2 \Big) \dd t
                   & - 2 \DS\int_0^1 \cW^{(1)}_t \cW^{(2)}_t \, \dd t \\
                   - 2 \DS\int_0^1 \cW^{(1)}_t \cW^{(2)}_t \, \dd t
                   & \DS\int_0^1 \Big( \big( \cW^{(1)}_t \big)^2
                                     + \big( \cW^{(2)}_t \big)^2 \Big) \dd t
                  \end{bmatrix}  =: \tbA^{(3)} ,
 \end{align*}
 where \ $(\cW^{(1)}_t)_{t\in\RR_+}$ \ and  \ $(\cW^{(2)}_t)_{t\in\RR_+}$ \ are independent standard Wiener processes.
Next we derive an (asymptotic) expansion for \ $\bd_n$ \ (defined in \eqref{10def_dn}).
First we examine \ $\bd_{2n}$, $n\in\NN$.
\ We have \ $M_k = X_k - X_{k-2}-\mu_\vare = \vare_k - \mu_\vare$, $k\in\NN$, \ hence separating the expectations
 we get
 \begin{align*}
  \sum_{k=1}^{2n} M_k X_{k-1}
  &= \sum_{k=1}^n (\vare_{2k} - \mu_\vare) V_k
     + \sum_{k=1}^{n-1} (\vare_{2k+1} - \mu_\vare) U_k \\
  &= \sum_{k=1}^n (\vare_{2k} - \mu_\vare) \big(k\mu_\vare + V_k - \EE(V_k)\big)
     + \sum_{k=1}^{n-1}
        (\vare_{2k+1} - \mu_\vare) \big(k\mu_\vare + U_k - \EE(U_k)\big) \\
  &= \mu_\vare \sum_{k=1}^n k (\vare_{2k} - \mu_\vare)
     + \mu_\vare \sum_{k=1}^{n-1} k (\vare_{2k+1} - \mu_\vare) \\
  &\phantom{\quad}
     + \sum_{k=1}^n (V_k - \EE(V_k)) (\vare_{2k} - \mu_\vare)
     + \sum_{k=1}^{n-1} (U_k - \EE(U_k)) (\vare_{2k+1} - \mu_\vare) \\
  &= d_{2n,1}^{(1)} (2n)^{3/2} + d_{2n,2}^{(1)} 2n ,
 \end{align*}
 where
 \begin{align*}
  d_{2n,1}^{(1)}
  &:= \frac{\mu_\vare}{(2n)^{3/2}} \sum_{k=1}^n k (\vare_{2k} - \mu_\vare)
      + \frac{\mu_\vare}{(2n)^{3/2}} \sum_{k=1}^{n-1} k (\vare_{2k+1} - \mu_\vare) \\
  &\distr \frac{\mu_\vare \sigma_\vare}{2^{3/2}}
          \int_0^1 t \, \dd\big( \cW_t^{(1)} + \cW_t^{(2)} \big)
   =: d_1^{(1)} , \\
  d_{2n,2}^{(1)}
  &:= \frac{1}{2n} \sum_{k=1}^n (V_k - \EE(V_k)) (\vare_{2k} - \mu_\vare)
      + \frac{1}{2n} \sum_{k=1}^{n-1} (U_k - \EE(U_k)) (\vare_{2k+1} - \mu_\vare) \\
  &\distr \frac{\sigma_\vare^2}{2} \int_0^1 \cW_t^{(1)} \, \dd\cW_t^{(2)}
          + \frac{\sigma_\vare^2}{2} \int_0^1 \cW_t^{(2)} \, \dd\cW_t^{(1)}
   =: d_2^{(1)} .
 \end{align*}
Indeed, by \eqref{01seged1}, \eqref{01seged_main}, the continuous mapping theorem
 (see, e.g., van der Vaart \cite[Theorem 2.3]{Vaa}), and It\^o's formula,
 \begin{align*}
  \frac{1}{n^{3/2}} \sum_{k=1}^n k (\vare_{2k} - \mu_\vare)
    & = \frac{1}{n^{3/2}} \sum_{k=1}^n \sum_{j=1}^k (\vare_{2k} - \mu_\vare)
      = \frac{1}{n^{3/2}} \sum_{j=1}^n \sum_{k=j}^n (\vare_{2k} - \mu_\vare) \\
    & = \frac{1}{n^{3/2}}
     \sum_{j=1}^n
      \left( \sum_{k=1}^n (\vare_{2k} - \mu_\vare)
             - \sum_{k=1}^{j-1}(\vare_{2k} - \mu_\vare) \right) \\
    & = \frac{1}{n^{3/2}}
      \sum_{j=1}^n \Big( \big(U_n - \EE(U_n)\big) -\big(U_{j-1} - \EE(U_{j-1})\big) \\
    &  = \frac{1}{n^{1/2}}(U_n-\EE(U_n))
        - \frac{1}{n^{3/2}}\sum_{j=1}^n (U_{j-1} - \EE(U_{j-1})) \\
    &\distr \sigma_\vare\cW_1^{(1)} - \sigma_\vare\int_0^1 \cW_s^{(1)} \, \dd s
              \ase \sigma_\vare \int_0^1 s \, \dd \cW_s^{(1)} ,
 \end{align*}
 and
 \begin{align}\label{01seged_D2}
  \begin{split}
   \frac{1}{n^{3/2}} & \sum_{k=1}^{n-1} k (\vare_{2k+1} - \mu_\vare)
        = \frac{1}{n^{3/2}} \sum_{j=1}^{n-1}
             \left( \sum_{k=1}^{n-1} (\vare_{2k+1} - \mu_\vare)
             - \sum_{k=1}^{j-1}(\vare_{2k+1} - \mu_\vare) \right)\\
       & = \frac{1}{n^{3/2}}  \sum_{j=1}^{n-1}
            \Big( (V_n-\EE(V_n))  - (\vare_1 - \mu_\vare)
                     - (V_{j-1}-\EE(V_{j-1})) + (\vare_1-\mu_\vare) \Big)\\
       & = \frac{n-1}{n^{3/2}}(V_n-\EE(V_n))
          - \frac{1}{n^{3/2}}\sum_{j=1}^{n-1} \big(V_{j-1} - \EE(V_{j-1})\big)  \\
       & \distr
           \sigma_\vare\cW_1^{(2)} - \sigma_\vare\int_0^1 \cW_s^{(2)} \, \dd s
           \ase \sigma_\vare \int_0^1 s \, \dd \cW_s^{(2)},
     \end{split}
 \end{align}
 where \ $\ase$ \ denotes equality almost surely.
Further,
 \begin{align*}
  d_{2n,2}^{(1)}
  &= \frac{1}{2n}
     \sum_{k=1}^n \sum_{j=1}^k (\vare_{2j-1} - \mu_\vare) (\vare_{2k} - \mu_\vare)
     + \frac{1}{2n}
       \sum_{k=1}^{n-1} \sum_{j=1}^k
        (\vare_{2j} - \mu_\vare) (\vare_{2k+1} - \mu_\vare) \\
  &= \frac{1}{2n}
     \sum_{i=1}^n (\vare_{2i} - \mu_\vare)
     \sum_{j=1}^n (\vare_{2j-1} - \mu_\vare)
  = \frac{1}{2}\frac{1}{\sqrt{n}}(U_n-\EE(U_n))\frac{1}{\sqrt{n}}(V_n-\EE(V_n)) \\
  & \distr \frac{\sigma_\vare^2}{2} \cW_1^{(1)} \cW_1^{(2)}
    = d_2^{(1)} .
 \end{align*}
In a similar way,
 \begin{align*}
  \sum_{k=1}^{2n} M_k X_{k-2}
  &= \sum_{k=1}^{n-1} (\vare_{2k+1} - \mu_\vare) V_k
     + \sum_{k=1}^{n-1} (\vare_{2k+2} - \mu_\vare) U_k \\
  &= \sum_{k=1}^{n-1} (\vare_{2k+1} - \mu_\vare) \big(k\mu_\vare + V_k - \EE(V_k)\big)
     + \sum_{k=1}^{n-1}
        (\vare_{2k+2} - \mu_\vare) \big(k\mu_\vare + U_k - \EE(U_k)\big) \\
  &= \mu_\vare \sum_{k=1}^{n-1} k (\vare_{2k+1} - \mu_\vare)
     + \mu_\vare \sum_{k=1}^{n-1} k (\vare_{2k+2} - \mu_\vare) \\
  &\phantom{\quad}
     + \sum_{k=1}^{n-1} (V_k - \EE(V_k)) (\vare_{2k+1} - \mu_\vare)
     + \sum_{k=1}^{n-1} (U_k - \EE(U_k)) (\vare_{2k+2} - \mu_\vare) \\
  &= d_{2n,1}^{(2)} (2n)^{3/2} + d_{2n,2}^{(2)}2n ,
 \end{align*}
 where
 \begin{align*}
  d_{2n,1}^{(2)}
  &:= \frac{\mu_\vare}{(2n)^{3/2}} \sum_{k=1}^{n-1} k (\vare_{2k+1} - \mu_\vare)
      + \frac{\mu_\vare}{(2n)^{3/2}} \sum_{k=1}^{n-1} k (\vare_{2k+2} - \mu_\vare) \\
  &\distr \frac{\mu_\vare \sigma_\vare}{2^{3/2}}
          \int_0^1 t \, \dd\big( \cW_t^{(1)} + \cW_t^{(2)} \big)
   =: d_1^{(2)} = d_1^{(1)}, \\
  d_{2n,2}^{(2)}
  &:= \frac{1}{2n} \sum_{k=1}^{n-1} (V_k - \EE(V_k)) (\vare_{2k+1} - \mu_\vare)
      + \frac{1}{2n}
        \sum_{k=1}^{n-1} (U_k - \EE(U_k)) (\vare_{2k+2} - \mu_\vare) \\
  &\distr \frac{\sigma_\vare^2}{2} \int_0^1 \cW_t^{(1)} \, \dd\cW_t^{(1)}
          + \frac{\sigma_\vare^2}{2} \int_0^1 \cW_t^{(2)} \, \dd\cW_t^{(2)}
   =: d_2^{(2)} .
 \end{align*}
Indeed, by \eqref{01seged_D2}, we have
 \[
   \frac{1}{n^{3/2}}  \sum_{k=1}^{n-1} k (\vare_{2k+1} - \mu_\vare)
      \distr \sigma_\vare \int_0^1 s \, \dd \cW_s^{(2)},
 \]
 and using that
 \begin{align*}
   \frac{1}{n^{3/2}}  \sum_{k=1}^{n-1} k (\vare_{2k+2} - \mu_\vare)
        & =\frac{1}{n^{3/2}} \sum_{k=1}^{n-1} \sum_{j=1}^k (\vare_{2k+2} - \mu_\vare)
          =\frac{1}{n^{3/2}} \sum_{j=1}^{n-1} \sum_{k=j}^{n-1} (\vare_{2k+2} - \mu_\vare)  \\
        & = \frac{1}{n^{3/2}} \sum_{j=1}^{n-1}
             \left( \sum_{k=1}^{n-1} (\vare_{2k+2} - \mu_\vare)
             - \sum_{k=1}^{j-1}(\vare_{2k+2} - \mu_\vare) \right)\\
        & = \frac{1}{n^{3/2}} \sum_{j=1}^{n-1}
             \left( \sum_{k=2}^{n} (\vare_{2k} - \mu_\vare)
             - \sum_{k=2}^{j}(\vare_{2k} - \mu_\vare) \right),
 \end{align*}
 by \eqref{01seged_main}, the continuous mapping theorem and It\^o's formula, we get
 \begin{align}\nonumber
%  \begin{split}
    \frac{1}{n^{3/2}}  \sum_{k=1}^{n-1} k (\vare_{2k+2} - \mu_\vare)
       & = \frac{1}{n^{3/2}}  \sum_{j=1}^{n-1}
            \Big( (U_n-\EE(U_n))  - (\vare_2 - \mu_\vare)
                     - (U_{j}-\EE(U_{j})) + (\vare_2-\mu_\vare) \Big)\\ \label{01seged14}
       & = \frac{n-1}{n^{3/2}}(U_n-\EE(U_n))
          - \frac{1}{n^{3/2}}\sum_{j=1}^{n-1} \big(U_{j} - \EE(U_{j})\big)  \\ \nonumber
       & \distr
           \sigma_\vare\cW_1^{(1)} - \sigma_\vare\int_0^1 \cW_s^{(1)} \, \dd s
           \ase \sigma_\vare \int_0^1 s \, \dd \cW_s^{(1)}.
% \end{split}
 \end{align}
By similar arguments, using also the strong law of large numbers, we have
 \begin{align*}
  &\frac{1}{2n}
   \sum_{k=1}^{n-1} (V_k - \EE(V_k)) (\vare_{2k+1} - \mu_\vare)
   = \frac{1}{2n}
     \sum_{k=1}^{n-1} \sum_{j=1}^k
      (\vare_{2j-1} - \mu_\vare) (\vare_{2k+1} - \mu_\vare) \\
  &= \frac{1}{2n}
     \sum_{1 \leq j < k \leq n} (\vare_{2j-1} - \mu_\vare) (\vare_{2k-1} - \mu_\vare)
   = \frac{1}{4n}
     \left[ \left( \sum_{k=1}^n (\vare_{2k-1} - \mu_\vare) \right)^2
            - \sum_{k=1}^n (\vare_{2k-1} - \mu_\vare)^2 \right]\\
  & = \frac{1}{4}\left[\left(\frac{V_n-\EE(V_n)}{\sqrt{n}}\right)^2
                      - \frac{1}{n}\sum_{k=1}^n (\vare_{2k-1}-\mu_\vare)^2
                      \right] \\
  &\distr \frac{\sigma_\vare^2}{4} \left[ \big( \cW_1^{(2)} \big)^2 - 1 \right]
   \ase \frac{\sigma_\vare^2}{2} \int_0^1 \cW_t^{(2)} \, \dd\cW_t^{(2)} ,
 \end{align*}
 and
 \begin{align*}
  &\frac{1}{2n}
   \sum_{k=1}^{n-1} (U_k - \EE(U_k)) (\vare_{2k+2} - \mu_\vare)
   = \frac{1}{2n}
     \sum_{k=1}^{n-1} \sum_{j=1}^k
      (\vare_{2j} - \mu_\vare) (\vare_{2k+2} - \mu_\vare) \\
  &= \frac{1}{2n}
     \sum_{1 \leq j < k \leq n} (\vare_{2j} - \mu_\vare) (\vare_{2k} - \mu_\vare)
   = \frac{1}{4n}
     \left[ \left( \sum_{k=1}^n (\vare_{2k} - \mu_\vare) \right)^2
            - \sum_{k=1}^n (\vare_{2k} - \mu_\vare)^2 \right]\\
  & = \frac{1}{4}\left[\left(\frac{U_n-\EE(U_n)}{\sqrt{n}}\right)^2
                      - \frac{1}{n}\sum_{k=1}^n (\vare_{2k}-\mu_\vare)^2
                      \right] \\
  &\distr \frac{\sigma_\vare^2}{4} \left[ \big( \cW_1^{(1)} \big)^2 - 1 \right]
   \ase \frac{\sigma_\vare^2}{2} \int_0^1 \cW_t^{(1)} \, \dd\cW_t^{(1)} .
 \end{align*}

By Lemma \ref{Lem01_UV3}, \ $\bd_{2n,1}^{(1)}$, $\bd_{2n,2}^{(1)}$, $\bd_{2n,1}^{(2)}$ \
 and \ $\bd_{2n,2}^{(2)}$ \ also converge jointly in distribution as \ $n\to\infty$.
\ Hence we conclude
 \[
   \bd_{2n} = \bd_{2n,1} (2n)^{3/2} + \bd_{2n,2} 2n, \qquad n\in\NN,
 \]
 with
 \begin{align*}
  \bd_{2n,1}
     & := \begin{bmatrix}
         d_{2n,1}^{(1)} \\
         d_{2n,1}^{(2)} \\
       \end{bmatrix}
     = \frac{\mu_\vare}{(2n)^{3/2}}
      \begin{bmatrix}
       \sum_{k=1}^n k (\vare_{2k} - \mu_\vare)
       + \sum_{k=1}^{n-1} k (\vare_{2k+1} - \mu_\vare) \\[2mm]
       \sum_{k=1}^{n-1} k (\vare_{2k+1} - \mu_\vare)
       + \sum_{k=1}^{n-1} k (\vare_{2k+2} - \mu_\vare)
      \end{bmatrix} \\
  &\distr
      \begin{bmatrix}
         d_{1}^{(1)} \\
         d_{1}^{(2)} \\
       \end{bmatrix}
       = \frac{\mu_\vare \sigma_\vare}{2^{3/2}}
          \int_0^1 t \, \dd\big( \cW_t^{(1)} + \cW_t^{(2)} \big)
          \begin{bmatrix} 1 \\ 1 \end{bmatrix} =: \bd^{(1)} ,
 \end{align*}
 and
 \begin{align*}
  \bd_{2n,2}
  & := \begin{bmatrix}
         d_{2n,2}^{(1)} \\
         d_{2n,2}^{(2)} \\
       \end{bmatrix}
   = \frac{1}{2n}
      \begin{bmatrix}
        \sum_{k=1}^n (V_k - \EE(V_k)) (\vare_{2k} - \mu_\vare)
       + \sum_{k=1}^{n-1} (U_k - \EE(U_k)) (\vare_{2k+1} - \mu_\vare) \\[2mm]
        \sum_{k=1}^{n-1} (V_k - \EE(V_k)) (\vare_{2k+1} - \mu_\vare)
       +  \sum_{k=1}^{n-1} (U_k - \EE(U_k)) (\vare_{2k+2} - \mu_\vare)
      \end{bmatrix} \\
  &\distr
        \begin{bmatrix}
         d_{2}^{(1)} \\
         d_{2}^{(2)} \\
       \end{bmatrix}
        = \frac{\sigma_\vare^2}{2}
          \begin{bmatrix}
           \int_0^1 \cW_t^{(1)} \, \dd\cW_t^{(2)}
           + \int_0^1 \cW_t^{(2)} \, \dd\cW_t^{(1)}  \\[2mm]
           \int_0^1 \cW_t^{(1)} \, \dd\cW_t^{(1)}
           + \int_0^1 \cW_t^{(2)} \, \dd\cW_t^{(2)}
          \end{bmatrix}  =: \bd^{(2)} .
 \end{align*}
Similar expansion can be derived for \ $\bd_{2n+1}$, $n\in\ZZ_+$.
\ Namely,
 \[
   \bd_{2n+1} = \bd_{2n+1,1} (2n+1)^{3/2} + \bd_{2n+1,2}(2n+1),\qquad n\in\ZZ_+,
 \]
 with
 \begin{align*}
  \bd_{2n+1,1}
      := \begin{bmatrix}
         d_{2n+1,1}^{(1)} \\
         d_{2n+1,1}^{(2)} \\
       \end{bmatrix}
     = \frac{1}{(2n+1)^{3/2}}
      \begin{bmatrix}
        (2n)^{3/2}d_{2n,1}^{(1)} +  n\mu_\vare (\vare_{2n+1} - \mu_\vare) \\[2mm]
        (2n)^{3/2}d_{2n,1}^{(2)} +  n\mu_\vare (\vare_{2n+1} - \mu_\vare)
      \end{bmatrix}
  \distr
      \begin{bmatrix}
         d_{1}^{(1)} \\
         d_{1}^{(2)} \\
       \end{bmatrix}
        = \bd^{(1)} ,
 \end{align*}
 and
 \begin{align*}
  \bd_{2n+1,2}
   := \begin{bmatrix}
         d_{2n+1,2}^{(1)} \\
         d_{2n+1,2}^{(2)} \\
       \end{bmatrix}
   = \frac{1}{2n+1}
      \begin{bmatrix}
         (2n)d_{2n,2}^{(1)} + (U_n - \EE(U_n)) (\vare_{2n+1} - \mu_\vare) \\[2mm]
         (2n)d_{2n,2}^{(2)} + (V_n - \EE(V_n)) (\vare_{2n+1} - \mu_\vare)
      \end{bmatrix}
  \distr
        \begin{bmatrix}
         d_{2}^{(1)} \\
         d_{2}^{(2)} \\
       \end{bmatrix}
        = \bd^{(2)} .
 \end{align*}
Indeed, \ $(\vare_{2n+1}-\mu_\vare)/\sqrt{n}\stoch 0$ \ as \ $n\to\infty$, \ and, by the independence
 of \ $U_n - \EE(U_n)$ \  and \ $\vare_{2n+1}-\mu_\vare$, \ we have
  \begin{align*}
   \EE\left(\frac{(U_n - \EE(U_n)) (\vare_{2n+1} - \mu_\vare)}{2n+1}\right)^2
      = \frac{\EE(U_n - \EE(U_n))^2 \EE(\vare_{2n+1} - \mu_\vare)^2}{(2n+1)^2}
      = \frac{n\sigma_\vare^4}{(2n+1)^2}
      \to 0
  \end{align*}
 as \ $n\to\infty$, \ which yields that \ $(U_n - \EE(U_n)) (\vare_{2n+1} - \mu_\vare)/(2n+1)\qmean 0$
 \ as \ $n\to\infty$, \ where \ $\qmean$ \ denotes convergence in \ $L_2$-sense.
Similarly, one can derive
  \ $(V_n - \EE(V_n)) (\vare_{2n+1} - \mu_\vare)/(2n+1)\qmean 0$ \ as \ $n\to\infty$.

Hence we have an (asymptotic) expansion for \ $\bd_{n}$, \ namely,
 \begin{align}\label{01seged3}
    \bd_{n} = \bd_{n,1} n^{3/2} + \bd_{n,2} n,
    \qquad n\in\NN,
 \end{align}
 where \ $\bd_{n,1}\distr \bd^{(1)}$ \ and \ $\bd_{n,2}\distr \bd^{(2)}$ \ as \ $n\to\infty$.

By \eqref{01seged2} and \eqref{01seged3},  we have the (asymptotic) expansion
 \begin{align}
  \begin{split}
   \tbA_n \bd_n
   &= ( \tbA_{n,1} n^3 + \tbA_{n,2} n^{5/2} + \tbA_{n,3} n^2 )
      ( \bd_{n,1} n^{3/2} + \bd_{n,2} n) \\
   &= \tbA_{n,1} \bd_{n,1} n^{9/2}
      + (\tbA_{n,1} \bd_{n,2} + \tbA_{n,2} \bd_{n,1}) n^4 \\
   &\phantom{\quad}
      + (\tbA_{n,2} \bd_{n,2} + \tbA_{n,3} \bd_{n,1}) n^{7/2}
      + \tbA_{n,3} \bd_{n,2} n^3 , \qquad n\in\NN.
  \end{split}
  \label{tAd}
 \end{align}
By Lemma \ref{Lem01_UV3}, Slutsky's lemma and the continuous mapping theorem,
 \ $\tbA_{n,i}$, $i=1,2,3$ \ and \ $\bd_{n,i}$, $i=1,2$ \ converge jointly in distribution,
 and hence the coefficients of the above expansion also converge jointly in distribution.
Further, we show that the first leading term has no influence by which we mean that
 \ $n^{1/2} \tbA_{n,1} \bd_{n,1} \stoch0$ \ as \ $n\to\infty$.
\ Indeed, we have
 \begin{align}\label{01seged17}
   &\tbA_{2n,1} =
               \frac{\mu_\vare^2}{12}
                 \begin{bmatrix} 1 & -1 \\ -1 & 1 \end{bmatrix}
               + \begin{bmatrix} O(n^{-1}) & O(n^{-1}) \\ O(n^{-1}) & O(n^{-1}) \end{bmatrix},
                  \qquad\text{$n\in\NN$,} \\[2mm]\label{01seged18}
   &\bd_{2n,1} = \frac{\mu_\vare}{(2n)^{3/2}}
                \sum_{k=1}^{n-1}
                  k ( \vare_{2k+1} - \mu_\vare + \vare_{2k+2}- \mu_\vare )
                  \begin{bmatrix} 1 \\ 1 \end{bmatrix}
           + \frac{\mu_\vare}{(2n)^{3/2}} \sum_{k=0}^{n-1}(\vare_{2k+2} - \mu_\vare)
                   \begin{bmatrix} 1 \\ 0 \end{bmatrix},
 \end{align}
 and hence
 \begin{align*}
  (2n)^{1/2} \tbA_{2n,1} \bd_{2n,1}
     & = \frac{\mu_\vare^3}{24n} \sum_{k=0}^{n-1}(\vare_{2k+2} - \mu_\vare)
                      \begin{bmatrix} 1 \\ -1 \end{bmatrix}\\
     &\phantom{=\;}  + \sum_{k=1}^{n-1}
                    k ( \vare_{2k+1} - \mu_\vare + \vare_{2k+2}- \mu_\vare )
                    \begin{bmatrix} O(n^{-2}) \\ O(n^{-2}) \end{bmatrix}\\
     &\phantom{=\;}
                  + \sum_{k=0}^{n-1}(\vare_{2k+2} - \mu_\vare)
           \begin{bmatrix} O(n^{-2}) \\ O(n^{-2}) \end{bmatrix},
     \qquad \text{$n\in\NN$.}
 \end{align*}
The above formulas with \ $O(n^{-1})$ \ and \ $O(n^{-2})$ \ are meant to be
 entrywise and coordinatewise, respectively.
Further, for sequences \ $(\zeta_n)_{n\in\NN}$, \ $(\eta_n)_{n\in\NN}$ \
 of real-valued random variables and a sequence \ $(\theta_n)_{n\in\NN}$
 \ of real numbers such that \ $\theta_n\ne0$, $n\in\NN$, \ the notation
 \ $\zeta_n=\eta_n O(\theta_n)$, $n\in\NN$, \ means that there exists a sequence
 \ $(\kappa_n)_{n\in\NN}$ \ of real numbers such that
 \ $\zeta_n=\eta_n \kappa_n$, $n\in\NN$, \ and
 \ $\sup_{n\in\NN}\vert\frac{\kappa_n}{\theta_n}\vert<\infty$.
\ By the strong law of large numbers,
 \begin{align*}
   \frac{1}{n} \sum_{k=0}^{n-1}(\vare_{2k+2} - \mu_\vare)
       \as \EE(\vare_2-\mu_\vare) = 0,
 \end{align*}
 and by \eqref{01seged_D2}, \eqref{01seged14} and Slutsky's lemma,
 \begin{align*}
   \frac{1}{n^2} \sum_{k=1}^{n-1}k(\vare_{2k+1} - \mu_\vare + \vare_{2k+2}- \mu_\vare)
      \stoch 0
   \qquad \text{as \ $n\to\infty$.}
 \end{align*}
Hence \ $(2n)^{1/2} \tbA_{2n,1} \bd_{2n,1} \stoch {\bf 0}\in\RR^2$.
\ Similarly one can prove that \ $\tbA_{2n+1,1} \bd_{2n+1,1}\stoch 0$ \
 as \ $n\to\infty$. \
Indeed, for all \ $n\in\ZZ_+$, \
 \begin{align*}
   &\tbA_{2n+1,1}\bd_{2n+1,1}
      = \begin{bmatrix}
          H_{2n+1,1} & -G_{2n+1,1} \\
          -G_{2n+1,1} & F_{2n+1,1} \\
        \end{bmatrix}
        \begin{bmatrix}
          d_{2n+1,1}^{(1)} \\
          d_{2n+1,1}^{(2)} \\
        \end{bmatrix}\\
   &\qquad =
        \frac{(2n)^3}{(2n+1)^{9/2}}
        \begin{bmatrix}
          H_{2n,1} + \frac{\mu_\vare^2}{8n} & -\left(G_{2n,1}+\frac{\mu_\vare^2}{8n}\right)\\
          -\left(G_{2n,1}+ \frac{\mu_\vare^2}{8n} \right) & F_{2n,1} + \frac{\mu_\vare^2}{8n}\\
        \end{bmatrix}
        \begin{bmatrix}
          (2n)^{3/2}d_{2n,1}^{(1)} + n\mu_\vare(\vare_{2n+1}-\mu_\vare) \\
          (2n)^{3/2}d_{2n,1}^{(2)} + n\mu_\vare(\vare_{2n+1}-\mu_\vare) \\
        \end{bmatrix}\\
   &\qquad = \frac{(2n)^3}{(2n+1)^{9/2}}
       \left(\tbA_{2n,1}
               + \frac{\mu_\vare^2}{8n}
               \begin{bmatrix}
                 1  & -1 \\
                 -1 & 1 \\
              \end{bmatrix}
        \right)
        \left((2n)^{3/2}\bd_{2n,1}
         + n\mu_\vare(\vare_{2n+1}-\mu_\vare)\begin{bmatrix}
                                               1 \\
                                               1 \\
                                             \end{bmatrix}
          \right),%\\
  \end{align*}
 and hence
 \begin{align*}
   \tbA_{2n+1,1}\bd_{2n+1,1}
    & = \left(\frac{2n}{2n+1}\right)^{9/2}  \tbA_{2n,1}\bd_{2n,1}
        + \frac{8n^4}{(2n+1)^{9/2}}\mu_\vare(\vare_{2n+1}-\mu_\vare)
           \tbA_{2n,1}  \begin{bmatrix}
                             1 \\
                             1 \\
                        \end{bmatrix} \\
   &\phantom{=\;}
        + \left(\frac{2n}{2n+1}\right)^{9/2}
         \frac{\mu_\vare^2}{8n}
             \begin{bmatrix}
                 1  & -1 \\
                 -1 & 1 \\
              \end{bmatrix}
           \bd_{2n,1},
 \end{align*}
 where, by \eqref{01seged17} and \eqref{01seged18},
 \begin{align*}
   &\tbA_{2n,1}
     \begin{bmatrix}
       1 \\
       1 \\
    \end{bmatrix}
    = O(n^{-1})
     \begin{bmatrix}
       1 \\
       1 \\
     \end{bmatrix} \qquad \text{as \ $n\to\infty$,}\\
  &  \begin{bmatrix}
                 1  & -1 \\
                 -1 & 1 \\
      \end{bmatrix}
           \bd_{2n,1}
     = \frac{\mu_\vare}{(2n)^{3/2}}
       \sum_{k=0}^{n-1}(\vare_{2k+2} - \mu_\vare)
        \begin{bmatrix}
       1 \\
       -1 \\
     \end{bmatrix}
     = \frac{\mu_\vare}{(2n)^{3/2}}
        (U_k - \EE(U_k))
       \begin{bmatrix}
       1 \\
       -1 \\
     \end{bmatrix}.
 \end{align*}
Using Lemma \ref{Lem01_UV3}, Slutsky's lemma and that \ $(2n)^{1/2}\tbA_{2n,1}\bd_{2n,1}\stoch {\bf 0}\in\RR^2$
 \ (which was proved earlier), we get \ $(2n+1)^{1/2}\tbA_{2n+1,1}\bd_{2n+1,1}\stoch {\bf 0}\in\RR^2$.
\ Then
 \[
     n^{1/2} \tbA_{n,1} \bd_{n,1} \stoch \begin{bmatrix}
                                           0 \\
                                           0 \\
                                         \end{bmatrix}
                     \qquad \text{as \ $n\to\infty$.}
 \]
Hence using also that the coefficients of the expansion of \ $\tbA_n \bd_n$ \
 converge jointly in distribution we obtain
 \[
   n^{-4} \tbA_n \bd_n \distr \tbA^{(1)} \bd^{(2)} + \tbA^{(2)} \bd^{(1)} .
 \]
Here \ $\tbA^{(2)} \bd^{(1)} = {\bf 0}\in\RR^2$ \ and
 \[
   \tbA^{(1)} \bd^{(2)}
   = \frac{\mu_\vare^2 \sigma_\vare^2}{24}
     \int_0^1 \big( \cW_t^{(1)} - \cW_t^{(2)} \big)
               \dd\big( \cW_t^{(1)} - \cW_t^{(2)} \big)
     \begin{bmatrix} -1 \\ 1 \end{bmatrix} .
 \]
In fact, by Lemma \ref{Lem01_UV3}, Slutsky's lemma and the continuous mapping theorem,
 we have joint convergence of \ $n^{-5}\det(\bA_n)$ \ and \ $n^{-4} \tbA_n \bd_n$, \ and hence,
 using also \eqref{01seged16}, we get

 \begin{align*}
  \big( n^{-5}& \det(\bA_n) , \, n^{-4} \tbA_n \bd_n \big) \\
  &\distr
  \left( \frac{\mu_\vare^2 \sigma_\vare^2}{24}
          \int_0^1 \big( \cW_t^{(1)} - \cW_t^{(2)} \big)^2 \dd t , \,
          \frac{\mu_\vare^2 \sigma_\vare^2}{24}
          \int_0^1 \big( \cW_t^{(1)} - \cW_t^{(2)} \big)
               \dd\big( \cW_t^{(1)} - \cW_t^{(2)} \big)
     \begin{bmatrix} -1 \\ 1 \end{bmatrix} \right) \\
   &\distre
   \left( \frac{\mu_\vare^2 \sigma_\vare^2}{12}
          \int_0^1 \big( \cW_t\big)^2 \dd t , \,
          \frac{\mu_\vare^2 \sigma_\vare^2}{12}
          \int_0^1  \cW_t \dd  \cW_t
     \begin{bmatrix} -1 \\ 1 \end{bmatrix} \right),
 \end{align*}
where \ $\distre$ \ means equality in distribution and the last step follows by that
 \ $2^{-1/2} \big( \cW_t^{(1)} - \cW_t^{(2)} \big)$, \ $t \in \RR_+$, \ is a standard Wiener process.

Let us consider the function \ $g$ \ defined in \eqref{01seged12}.
Since \ $g$ \ is continuous on \ $(\RR\setminus\{0\})\times\RR^2$ \ and
 \[
   \PP\left(
      \left( \frac{\mu_\vare^2\sigma_\vare^2}{12} \int_0^1 \big( \cW_t\big)^2 \dd t ,
           \frac{\mu_\vare^2\sigma_\vare^2}{12}\int_0^1  \cW_t \dd  \cW_t
              \begin{bmatrix} -1 \\ 1 \end{bmatrix}
    \right)  \in (\RR\setminus\{0\})\times\RR^2
   \right)=1,
 \]
 the continuous mapping theorem yields that
 \begin{align*}
   g\left(n^{-5}\det(\bA_n), n^{-4} \tbA_n \bd_n \right)
     \distr
    & g\left( \frac{\mu_\vare^2\sigma_\vare^2}{12} \int_0^1 \big( \cW_t\big)^2 \dd t ,
           \frac{\mu_\vare^2\sigma_\vare^2}{12}\int_0^1  \cW_t \dd  \cW_t
              \begin{bmatrix} -1 \\ 1 \end{bmatrix} \right)\\
    & \distre
       \frac{\int_0^1  \cW_t \dd  \cW_t}{\int_0^1 \big( \cW_t\big)^2 \dd t }
       \begin{bmatrix} -1 \\ 1 \end{bmatrix},
 \end{align*}
 where the last step follows by \ $\PP\left(\int_0^1 \big( \cW_t\big)^2 \dd t >0 \right)=1$.
\ By Proposition \ref{Pro1}, we have
 \begin{align*}
  \PP\left(
       n\begin{bmatrix}
                \halpha_n-\alpha \\
                \hbeta_n -\beta \\
              \end{bmatrix}
   = g\left(n^{-5}\det(\bA_n), n^{-4} \tbA_n \bd_n \right)
     \right)
  & \geq \PP\left(\sum_{k=1}^n X_{k-2}^2 > 0\right)\to1
   \qquad \text{as \ $n\to \infty$.}
 \end{align*}
Then Lemma \ref{Lem10_seged} concludes the proof.
\proofend

\end{document}